\documentclass[11pt,english]{article}
\usepackage[T1]{fontenc}
\pdfoutput=1
\usepackage[latin9]{inputenc}
\usepackage[a4paper]{geometry}
\usepackage{color}
\usepackage{babel}
\usepackage{verbatim}
\usepackage{float}
\usepackage{amsmath}
\usepackage{amsthm}
\usepackage{amssymb}
\usepackage{graphicx}
\usepackage{csquotes}
\usepackage{bm}
\usepackage{pgfplots}
\pgfplotsset{compat=newest}
\pgfplotsset{plot coordinates/math parser=false}
\usepackage[unicode=true,pdfusetitle,
 bookmarks=true,bookmarksnumbered=false,bookmarksopen=false,
 breaklinks=false,pdfborder={0 0 0},pdfborderstyle={},backref=false,colorlinks=false]
 {hyperref}

\makeatletter
\theoremstyle{definition}
\newtheorem*{example*}{\protect\examplename}
\theoremstyle{plain}
\newtheorem{theorem}{\protect\theoremname}[section] % remove [section] to make Theorems, Lemmas and Remarks be counted by section
\theoremstyle{plain}
\newtheorem{lemma}[theorem]{\protect\lemmaname}
\theoremstyle{plain}
\newtheorem{remark}[theorem]{\protect\remarkname}
\theoremstyle{plain}

\theoremstyle{plain}
\newtheorem{notation}[theorem]{\protect\notationname}

\numberwithin{equation}{section}

\makeatother

\providecommand{\examplename}{Example}
\providecommand{\lemmaname}{Lemma}
\providecommand{\theoremname}{Theorem}
\providecommand{\remarkname}{Remark}
\providecommand{\notationname}{Notation}
\providecommand{\propositionname}{Proposition}

\newcommand{\bh}{\mathbf{h}}
\newcommand{\bp}{\mathbf{p}}
\newcommand{\bg}{\mathbf{g}}

\newcommand{\dd}{\text{d}}
\newcommand{\RR}{\mathbb{R}}
\newcommand{\N}{\mathcal{N}}
\newcommand{\s}{\mathbb{S}}
\newcommand{\xx}{\mathbf{x}}
\newcommand{\yy}{{\mathbf{y}}}
\newcommand{\n}{\mathbf{n}}

\newcommand{\ee}{\mathbf{e}}
\newcommand{\kk}{\mathbf{k}}

\newcommand{\ka}{\kappa}

\newcommand{\Rr}{\mathcal{R}}
\newcommand{\Oo}{\mathcal{O}}
\newcommand{\F}{\mathbf{F}}
\newcommand{\zero}{\mathbf{0}}
\newcommand{\ii}{\text{i}}
\newcommand{\nn}{\mathbf{\bar n}}
\newcommand{\yp}{\mathbf{y}_{\text{p}}}
\newcommand{\Vector}[1]{{\left(\begin{matrix} #1 \end{matrix}\right)}}

\usepackage[
backend=biber,
style=numeric,
sorting=nyt,
giveninits=true
]{biblatex}
\begin{document}
\title{High order corrected trapezoidal rules for a class of singular integrals}
\author{Federico Izzo\footnote{Corresponding author}\ $^{,}$\footnote{Department of Mathematics, KTH Royal Institute of Technology, Stockholm, Sweden (\href{mailto:izzo@kth.se}{izzo@kth.se})}
\and Olof Runborg\footnote{Department of Mathematics, KTH Royal Institute of Technology, Stockholm, Sweden (\href{mailto:olofr@kth.se}{olofr@kth.se})}
\and Richard Tsai\footnote{ Department of Mathematics and Oden Institute for Computational Engineering and Sciences, The University of Texas at Austin, Austin TX, USA (\href{mailto:ytsai@math.utexas.edu}{ytsai@math.utexas.edu})}}
\date{}
\maketitle
\begin{abstract}
We present a family of high order trapezoidal rule-based quadratures for a class of singular integrals, where the integrand has a point singularity.
The singular part of the integrand is expanded in a Taylor series involving terms of increasing smoothness.
The quadratures are based on the trapezoidal rule, with the quadrature weights for Cartesian nodes close to the singularity judiciously corrected based on the expansion. High order accuracy can be achieved by utilizing a sufficient number of correction nodes around the singularity to approximate the terms in the series expansion.
The derived quadratures are applied to the Implicit Boundary Integral formulation of surface integrals involving the Laplace layer kernels.\\


\end{abstract}

\section{Introduction}

The trapezoidal rule is a simple and robust algorithm for approximating integrals.
{In general, it has
second order accuracy, 
but when applied to compactly supported smooth functions, the accuracy
is much higher.
However, 
as the integrand becomes less smooth,
the accuracy deteriorates,
which makes the method unsuitable for
singular integrals
such as those found in boundary integral equations.} 
This paper develops a systematic approach to derive high order corrected trapezoidal rules for integrals involving a class of integrands
that are singular at one point.

Let $f:\mathbb{R}^n\setminus\{\zero\}\mapsto \mathbb{R}$ be a compactly supported 
function with an integrable singularity at $\zero$. 
A crude way to approximate its integral is with 
the ``punctured'' trapezoidal rule $T_h^0$, 
where $h$ denotes the discretization parameter. It
equals the
standard trapezoidal rule, but sets $f\equiv 0$ in a
in a small $h$-dependent region $\N_h$ surrounding
the singular point
\[
T_h^{0}[f]:=h^n\sum_{\mathbf{y}\in h\mathbb{Z}^n \setminus \N_h}f(\mathbf{y}).
\]
This gives a low order accurate approximation.
With $\Rr_h[f]$ denoting the error in the quadrature rule, we write 
\begin{equation}\label{eq:error-Rh}
    \int f(\xx)\dd\xx = T_h^{0}[f]+\Rr_h[f]. 
\end{equation}
One direction to improve the accuracy is to 
add back the function values {at} the 
excluded points in $\N_h$ with
judiciously chosen weights, such that they
well approximate 
$\Rr_h[f]$. This can be seen as a
correction of the standard trapezoidal rule,
locally around the singularity. 
The overall simplicity of the method is therefore maintained.
The approach has been used for example in \cite{kapur1997high,marin2014corrected,wumartinsson2020zeta,wumartinsson2020corrected} for the trapezoidal rule and in \cite{strain1995locally,lether1991numerical} for other quadrature methods.

\begin{figure}
    \begin{center}
        \includegraphics[scale=1]{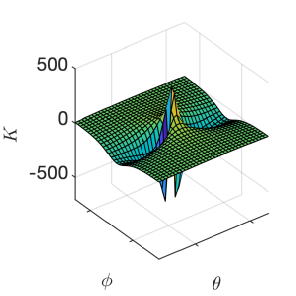}\hspace{0.5cm}\includegraphics[scale=1]{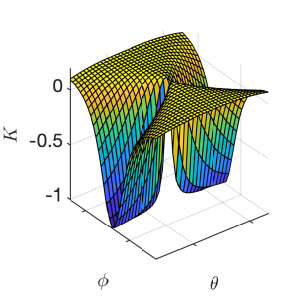}
    \end{center}
    \caption{\textbf{Singular behavior of a layer kernel}. Two plots related to the double-layer kernel $\frac{\partial G}{\partial \n_y}(\bar\xx, \bar\yy)=\frac{\partial G}{\partial \n_y}(\mathbf{v}(0,0),\mathbf{v}(\theta,\phi))$ from boundary integral formulations, where $\mathbf{v}(\theta,\phi)$ is a surface parametrization centered around $\bar\xx=\mathbf{v}(0,0)$. The kernel is in the form \eqref{eq:s-singular-integrand} using $(\theta,\phi)$ near its singular point.
    On the left, we plot the kernel $\frac{\partial G}{\partial \n_y}(\mathbf{v}(0,0),\mathbf{v}(\theta,\phi))$ on a uniform grid around $(0,0)$. On the right we plot  $\ell(\vert\mathbf{w}\vert,\mathbf{w}/\vert\mathbf{w}\vert)$, $\mathbf{w}:=(\theta,\phi)$, obtained by multiplying the same kernel by $\vert\mathbf{w}\vert$. We can clearly see that $\ell(0,\mathbf{u})$ is not constant in $\mathbf{u}$.}
    \label{fig:DL-example-plotting}
\end{figure}

In this article, we consider two-dimensional singular integrands
$f=s\,v$, where $v\in C^\infty_c(\RR^2)$ 
and $s$ is 
of the 
following form: 
\begin{equation}\label{eq:s-singular-integrand}
    s(\xx)=\dfrac{1}{\vert\xx\vert}\ell\left(\vert\xx\vert,\dfrac{\xx}{\vert\xx\vert}\right), 
\end{equation}
for some 
smooth function
$\ell:\RR\times\s^1\to\RR$.
Nevertheless, $\ell\left(\vert\xx\vert,\xx/\vert\xx\vert\right)$ is not necessarily a smooth function of $\xx$ in $\zero$ if $\ell(0,\mathbf{u})$ is non-constant in $\mathbf{u}$.
Singular functions of the type \eqref{eq:s-singular-integrand}
are found in many applications. For instance, 
if $\mathbf{g}:\RR^2\to\RR^3$ has 
a simple zero {at} the origin,
then $s(\xx)=1/\vert\mathbf{g(\xx)}\vert$ is of this type, as proven
in Lemma~\ref{lem:g=f-over-distance}. {They also
characterize the singular behavior of the kernels found in the boundary integral equations for elliptic problems.}
In three dimensions, the kernel is a function of two spatial variables, $\bar\xx,\bar{\mathbf{y}}$, but the integral typically involves the product of the kernel and a smooth function over a smooth and compact surface.
In this setup, the singularity in the integrand depends on $1/\vert\bar\xx-\bar{\mathbf{y}}\vert$ and also on ``the angle of approach'', which corresponds to the way $\bar{\mathbf{y}}$ approaches $\bar\xx$ along a two dimensional surface.
See Figure \ref{fig:DL-example-plotting} for an illustration of the singular behavior of boundary integrals layer kernels.
{Note that we limit ourselves to 
compactly supported integrands. These can also
be seen as the restrictions of periodic
functions which are smooth away from the point
singularity. If the integrand is not zero
at the boundary of the integration domain, 
additional boundary corrections must also be
introduced; see discussions in \cite{kapur1997high, aguilar2002high}.
}

In \cite{izzo2021corrected}, we derived a second order accurate method. In this paper, we generalize our approach systematically to derive higher order methods.
Our approach is to Taylor expand the function $\ell$
in its first argument,
and recognize that 
the smoothness of the remainder term increases
with order, and can eventually be integrated accurately 
with the standard
trapezoidal rule. We therefore only need to derive corrections
for the leading Taylor terms, which are all
of the form $\vert\xx\vert^j\phi(\xx/\vert\xx\vert)$,
for some $j$ and $\phi:\s^1\to\RR$. The details are presented in Section~\ref{sec:corrected-trapezoidal-rules}.

One of the main motivations for the proposed approach is to
provide the Implicit Boundary Integral Methods (IBIMs, see \cite{kublik2013implicit}) with high order convergent quadratures. IBIMs are volumetric integral formulations of 
classical boundary integrals and do not rely on explicit parameterization of surfaces (the ``boundary'' in the boundary integrals). {The IBIM approach gives a way to compute accurate surface integrals, integral equations and variational problems on surfaces for other non-parametric methods, including the level set methods, e.g. \cite{osher1988fronts,sethian1999level,osher2006level,chu2018volumetric} , and the closest point methods, e.g. \cite{ruuth2008simple,macdonald2008level}.}

In Section~\ref{sec:application-layer-potentials-3D}, we apply our new high order corrected trapezoidal rules to the singular integrals derived from IBIMs. In that formulation, the integrand is singular along a line and for each fixed plane it has a point singularity. To compute the volumetric integrals from IBIMs, our quadrature rules for integration in two dimensions are therefore applied plane by plane; see Section~\ref{sec:plane-by-plane}. We show in Theorem~\ref{thm:kernelexpansions} that the
resulting singularity on each plane is of the type in
\eqref{eq:s-singular-integrand}, and derive explicit expressions for the required terms in the expansion. 
Some efforts are needed to extract the needed geometrical information of the surface. Specifically, intrinsic information about the surface (principal directions and curvatures, and third derivatives of its local representation) together with extrinsic information (signed distance function to the surface) are needed to apply the quadrature rule in addition to the information needed for the IBIM formulation.

Finally, numerical simulations for selected problems in two and three dimensions are presented in Section~\ref{sec:numerical-tests}. 

\section{The corrected trapezoidal rules}
\label{sec:corrected-trapezoidal-rules}

The standard trapezoidal rule
has a low order accuracy when
applied to singular integrals. 
In this section we show how one can
raise the order of accuracy for
integrands that are singular at a point,
by
correcting the computations {at a} few
grid points close to the singularity.
This type of corrections have
been applied successfully in a few settings earlier. See for example \cite{kapur1997high, marin2014corrected, wumartinsson2020zeta,wumartinsson2020corrected}.

We begin by defining the trapezoidal rules that
we will work with. 
Let $f$ be an integrable, {compactly
supported,} function on $\mathbb{R}^n$.
We are interested in approximating the integral $\int_{{\mathbb{R}}^n}f(\xx)\text{d}\xx$ by summation of the values of $f$ on the uniform grid $h\mathbb{Z}^n$. {Since} $f$ is supported in a compact set, the standard \emph{trapezoidal rule} becomes the following simple Riemann sum:
\begin{equation}
  T_{h}[f] := h^n\sum_{\mathbf{y}\in h\mathbb{Z}^n} f(\mathbf{y}) . \label{eq:trapezoidalrule_nD}
\end{equation}
{
In this case the order of accuracy
of the approximation is only limited
by the regularity of $f$.
If $f\in C_c^p(\mathbb{R}^n)$,
the error is at worst  $\mathcal{O}(h^p)$. 
See
e.g. the discussion and proofs in \cite{rauchtaylor2010}.
In particular, the trapezoidal rule enjoys spectral accuracy if $f\in C_c^\infty(\mathbb{R}^n)$. 
Here $C_c^p(\mathbb{R}^n)$
denotes the space of 
compactly supported functions
on $\mathbb{R}^n$ whose partial derivatives up to
order $p$ are continuous
(of all orders, if $p=\infty$).
}

If $f$ is smooth in ${\mathbb{R}}^n\setminus \{\xx_0\}$, singular at $\xx_0$, and $\int_{{\mathbb{R}}^n}f(\xx)\text{d}\xx$ exists as a Cauchy principal value, it is natural to modify the trapezoidal rule by
excluding the summation over some grid nodes close to $\xx_0$. We define the {\it punctured trapezoidal rule} with respect to $\N_h$ as
\begin{equation}
T_{h,\,\mathcal{N}_h}^{0}[f]:=h^n\sum_{\mathbf{y}\in h\mathbb{Z}^n \setminus \N_h(\xx_0)}f(\mathbf{y}), \label{eq:puncturedtrapez}
\end{equation}
where $\N_h(\xx_0)$ defines a small neighborhood around $\xx_0$, the region being ``punctured'' from $\mathbb{R}^n$. When $\xx_0$ lies on a grid node, one typically sets $\N_h(\xx_0)=\{\xx_0\}$; i.e. only the singularity point is removed from the standard trapezoidal rule. If $\xx_0$ does not lie on a grid node, one option is to remove the grid node $\xx_h$ that is closest to it. In this case, $\N_h(\xx_0)=\{\xx_h\}$. In general, $\N_h(\xx_0)$ may contain several grid nodes, although the number is typically finite and independent of $h$. We will write $\mathcal{N}_{h,m}$ to indicate that the set contains $m$ nodes.

Here we consider an integrand
that is the
product of a smooth factor
$v$ and a singular factor $s$,
which takes the form \eqref{eq:s-singular-integrand} near the origin.
The punctured trapezoidal rule converges for 
such singular functions, albeit with a lower rate.
In the case $\N_h(\mathbf{x}_0)=\{\xx_h\}$, we have the following theorem.
\begin{theorem}\label{thm:punctured-tr-s-ell-around-singularity}
Suppose $v\in C_c^\infty({\mathbb{R}}^n)$
and 
$s(\,\cdot-\xx_0)\,v(\cdot)\in C^\infty_c({\mathbb{R}}^n\setminus\{\xx_0\})$ {for any $\xx_0\in{\mathbb{R}}^n$.}
Assume furthermore that
for some $r_0>0$ there exist $j\in\mathbb{Z}$ and $\ell\in C^\infty((-r_0,r_0)\times\s^{n-1})$ such that
\[
    s(\xx)=\vert\xx\vert^j\ell\left(\vert\xx\vert,\dfrac{\xx}{\vert\xx\vert}\right),\qquad \xx\in B_{r_0}(\zero).
\]
Then, for $j\geq 1-n$, 
\begin{equation*}
    \left\vert\int_{{\mathbb{R}}^n} s(\xx-\xx_0)v(\xx) \text{\emph{d}}\xx - T^0_{h,\,\N_h}[s(\,\cdot-\xx_0)\, v(\cdot)]\right\vert \leq C h^{j+n},
\end{equation*}
where the constant $C$ is independent of $h$, but depends on $j$, $\ell$ and $v$. 
\end{theorem}
\noindent The proof is given in Appendix~\ref{sec:A:thm-s-singularity}, where without loss of generality we consider $\xx_0=\zero$ and $\N_h(\zero)=\{\zero\}$.

We now give a brief summary of the steps that 
we shall take in Sections ~\ref{sub:corrections-splitting-intro-tr}--\ref{sub:additive-splitting} 
to correct the trapezoidal rule
for the two-dimensional case $n=2$ and $j=-1$ in Theorem~\ref{thm:punctured-tr-s-ell-around-singularity}. 

In Section~\ref{sub:corrections-splitting-intro-tr} 
we expand $\ell$ in its first
argument to derive {a series}  
of the form
\begin{align*}
s(\xx)\equiv&\, \dfrac{1}{\vert\xx\vert}\ell\left(\vert\xx\vert,\dfrac{\xx}{\vert\xx\vert}\right)  
=\sum_{k=0}^q s_k(\xx)  
+ \triangle_q s(\xx),
\end{align*}
for some 
functions $s_k$ and $\triangle_q s$.
Theorem~\ref{thm:punctured-tr-s-ell-around-singularity} states that 
the error $\Rr_h$, as defined in \eqref{eq:error-Rh}, in applying $T^0_{h,\,\mathcal{N}_h}$ to integrate $s_k$ is bounded above by $h^{k+1}$. More precisely,
$\Rr_h[s_k(\,\cdot-\xx_0)\,v(\,\cdot\,)]\sim \Oo(h^{k+1})$.

In Section~\ref{sub:first-order} we derive 
a weight $\omega$ for approximating the error $\mathcal{R}_h[s_k(\,\cdot-\xx_0)\,v(\,\cdot\,)]$. Multiplication of the weight by any smooth function $v$ should yield
\[
    \mathcal{R}_h[s_k(\,\cdot-\xx_0)\,v(\,\cdot\,)] = h^{k+1}\omega\,v(\xx_{h}) + \Oo(h^{k+2}),
\]
where $\xx_h$ is the grid node in $h\mathbb{Z}^2$ closest to $\xx_0$.
In addition, the weight depends on $s_k$ but not on $h$ and
$v$.
With this weight, we define
\[
Q_h^1[s_k(\,\cdot-\xx_0)\,v(\,\cdot\,)]:=\,T_{h,\,\mathcal{N}_{h,1}}^0 [s_k(\,\cdot-\xx_0)\,v(\,\cdot\,)] + h^{k+1}\omega\,v(\xx_{h}).
\]
Consequently, 
\begin{align*}
\int_{{\mathbb{R}}^2}s_k(\xx-\xx_0)v(\xx)\dd\xx =&\,  Q_h^1[s_k(\,\cdot-\xx_0)\,v(\,\cdot\,)] + \Oo(h^{k+2}).
\end{align*}

We then generalize this approach systematically in Sections~\ref{sub:second-order-correction}~and~\ref{sub:high-order-correction}. 
Eventually, we obtain quadratures 
\[
Q_h^p[s_k(\,\cdot-\xx_0)\,v(\,\cdot\,)] 
:=\,T_{h,\,\mathcal{N}_{h,\tilde p}}^0 [s_k(\,\cdot-\xx_0)\,v(\,\cdot\,)] + h^{k+1}\sum_{i=1}^{\tilde p}\omega_i\,v(\xx_{h,i}),
\]
such that
\[
\int_{{\mathbb{R}}^2}s_k(\xx-\xx_0)v(\xx)\dd\xx =\, Q_h^p[s_k(\,\cdot-\xx_0)\,v(\,\cdot\,)] + \Oo(h^{k+1+p}),
\]
formally for any $k\ge 0$. Here $\tilde p\ge p$ is a constant 
and $\xx_{h,i}$ are grid points near $\xx_0$.
These will be described more carefully later.  

\begin{remark}\label{rem:Q1-different-for-k}
The quadrature rule $Q_h^p[s_k]$ depends on the value of $k$ in the subscript of $s_k$, in addition to the function $s_k$, but for simplicity of notation we will not make this distinction.
\end{remark}

Finally, in Section~\ref{sub:additive-splitting} we 
combine the quadratures $ Q^p_h$ for $s_k\,v$ to define a quadrature~$\mathcal{U}^{p}_h$ of order $p\geq 2$ for the function $s\,v$ (recall that $s$ is expanded into a sum of $s_k$ for
$k=0,1,\cdots, q,$ and $\triangle_q s$). 
The order $p$ specifies how many expansions terms are needed ($q=p-2$) and which quadratures derived from correcting the punctured trapezoidal rules are needed for each term ($Q^{p-1-k}_h$ for $s_k$, and $T_{h,\,\mathcal{N}_{h,1}}^0$ for $\triangle_{q}s=\triangle_{p-2}s$):
\begin{align*}
    \int_{{\mathbb{R}}^2}s(\xx-\xx_0)v(\xx)\dd\xx =&\ \mathcal{U}^p_h[s(\,\cdot-\xx_0)v(\,\cdot\,)] + \Oo(h^p)\\
    =&\ \sum_{k=0}^{p-2} Q_h^{p-k-1}[s_k(\,\cdot-\xx_0)\,v(\,\cdot\,)] \\ &\ + T_{h,\,\mathcal{N}_{h,1}}^0[\triangle_{p-2}s(\,\cdot-\xx_0)\,v(\,\cdot\,)] + \Oo(h^p).
\end{align*}

\subsection{Expansion of the singular function}
\label{sub:corrections-splitting-intro-tr}

To integrate the function $s(\xx)=\vert\xx\vert^{-1}\ell(\vert\xx\vert,\xx/\vert\xx\vert)$ in \eqref{eq:s-singular-integrand} with high order accuracy, we use a \emph{divide et impera} strategy. For any $\mathbf{u}\in\s^1$, we expand $\ell(r,\mathbf{u})$ with respect to the first variable, and approach each of the expansion components separately: $\ell$ becomes 
$$
\ell(r,\mathbf{u}) = \ell(0,\mathbf{u})+r\,\partial_r\ell(0,\mathbf{u})+\dfrac{1}{2}r^2\partial_r^2\ell(0,\mathbf{u})
+\cdots,
$$
and we write $s$ formally as the series
\begin{align}
    &s(\xx) = s_0(\xx) + s_1(\xx) + s_2(\xx) + \cdots \nonumber\\[0.2cm]
    & \hspace{0.77cm} = \dfrac{1}{\vert\xx\vert}\,\phi_0\left(\dfrac{\xx}{\vert\xx\vert}\right)+\phi_1\left(\dfrac{\xx}{\vert\xx\vert}\right) + \vert\xx\vert\,\phi_2\left(\dfrac{\xx}{\vert\xx\vert}\right) + \cdots\label{eq:s-expansion}\\[0.3cm]
    &\text{where }\ s_k(\xx):=\vert\xx\vert^{k-1}\phi_k\left(\frac{\xx}{\vert\xx\vert}\right), \text{ and }\  \phi_k\left(\dfrac{\xx}{\vert\xx\vert}\right):=\dfrac{1}{k!}\partial_r^k\ell\left(0,\dfrac{\xx}{\vert\xx\vert}\right)\nonumber.
\end{align}

If we use $q$ terms in this expansion, 
we expect that these terms $s_k$ strip away the singularity in $s$ at $\xx=\zero$ so that what is left behind from the expansion, i.e. the remainder term 
\begin{equation}\label{eq:ell-remainder}
\triangle_q s(\xx)\,:=\,s(\xx)-(s_0(\xx)+s_1(\xx)+\cdots+s_q(\xx))
\end{equation}
can be approximated directly with the (unmodified) trapezoidal rule and achieve the order of accuracy desired without needing special quadrature.

This property is expressed in the following lemma.
\begin{lemma}\label{lem:f-over-distance-expansion}
Let $s$ be of the kind \eqref{eq:s-singular-integrand}. Let $r_0>0$ be such that $\ell\in C^\infty((-r_0,r_0)\times\s^1)$. For any integer $q\ge 0$,
there exist $\sigma:{\mathbb{R}}\times\s^1\to{\mathbb{R}}$ such that $\sigma \in C^\infty((-r_0,r_0)\times \s^1)$ and
\begin{equation}\label{eq:R-remainder-form}
    \triangle_q s(\xx)=\vert\xx\vert^{q}{\sigma}(\vert\xx\vert,\xx/\vert\xx\vert).
\end{equation}
\end{lemma}
The proof of this lemma can be found in Appendix~\ref{sub:A:f-over-distance-expansion}. {From this result and the previous Theorem~\ref{thm:punctured-tr-s-ell-around-singularity} we can express the following lemma.}
\begin{lemma}\label{lem:Rq-remainder-trapezoidal-rule}
The term $\triangle_q s$ in \eqref{eq:ell-remainder} is integrated by the punctured trapezoidal rule~\eqref{eq:puncturedtrapez} with order $q+2$:
\begin{equation*}
    \left\vert\int_{{\mathbb{R}}^2} \triangle_q s(\xx-\xx_0)v(\xx) \dd\xx - T^0_{h,\,\mathcal{N}_h} [\triangle_q s(\cdot-\xx_0)v(\,\cdot\,)]\right\vert \leq C h^{q+2}. 
\end{equation*}
\end{lemma}
Hence, to get high order, it is sufficient
to derive corrected trapezoidal rules for the expansion terms 
\begin{equation}\label{eq:sk-integrand}
  s_k(\xx)=\vert\xx\vert^{k-1}\,\phi_k\left(\dfrac{\xx}{\vert\xx\vert}\right), \qquad k=0,1,2,3,\dots.
\end{equation}
We start from first order correction for $s_k$ in Section~\ref{sub:first-order} and end with a general description for an arbitrarily high order method for $s$ in Section~\ref{sub:additive-splitting}.

%%%%%%%%%%%%%%%%%%%%%%%%%%%%%%%
\subsection{First order correction}
\label{sub:first-order}
Our goal in this section is to derive a first order in $h$ correction for the punctured trapezoidal rule applied to
\begin{equation}\label{thirdordererror}
  \int_{{\mathbb{R}}^2} s_k(\xx-\xx_0){v(\xx)}\dd\mathbf{x},
\end{equation}
for $s_k$ of the form \eqref{eq:sk-integrand}.
We assume $\phi_k\in C^\infty({\mathbb S^1})$, $k\geq 0$, and $v\in C^\infty_c({\mathbb R}^2)$. This correction will yield an error with its largest part proportional to $h^{k+2}$.

\subsubsection{The singular point rests on a grid node}
Without loss of generality, we assume that $\xx_0=\zero$ and lies on a grid node.  For such cases, the set $\N_h$ typically contains only the grid node where the singularity is. In our case, $\N_h(\zero)=\{\zero\}$. 
The smoothness of $s_k$ 
in \eqref{eq:sk-integrand}
increases with $k$.
Theorem~\ref{thm:punctured-tr-s-ell-around-singularity} tells us that for a function of this kind (two-dimensional, $j=k-1$) the error behaves as:
\[
\left\vert\int_{{\mathbb{R}}^2}s_k(\xx)v(\xx)\dd\xx - T_{h,\,\N_{h}}[s_k\,v] \right\vert \leq Ch^{k+1}.
\]
Following \cite{marin2014corrected}, one can show that the error has the form
\begin{equation*}
  \int_{{\mathbb{R}}^2}s_k(\xx){v(\xx)}\dd\xx =T_{h,\,\mathcal{N}_h}^0[s_k\, v] + h^{k+1}\,\omega[s_k] v(\mathbf{\mathbf{0}}) +\mathcal{O}(h^{k+2}),
\end{equation*}
where $\omega[s_k]$ is a constant independent of $v$ and $h$. 
In \cite{marin2014corrected} this is proven for $s_0(\xx)=1/\vert\xx\vert$ ($k=0$ and $\phi_0\equiv 1$). 

Hence we define the \emph{first order correction} $Q_h^1$ to the punctured trapezoidal rule as
\begin{equation}\label{eq:Q1-symmetric-grid-Marin}
Q_h^1[s_k\, v]:=T_{h,\,\mathcal{N}_h}^0[s_k\,v] + h^{k+1}\,\omega\left[s_k\right]v(\mathbf{0}).
\end{equation}
This quadrature rule thus
corrects the trapezoidal rule in
one node, the origin. It
will then have
an error of size $\Oo(h^{k+2})$.
\begin{remark}
In the special case when
$\phi_k\equiv 1$, due to symmetry with respect to the grid node at $\zero$,
the $\Oo(h^{k+2})$ terms cancel out, and
$Q_h^1$ achieves an accuracy of
$\Oo(h^{k+3})$.
\end{remark}

To find the weight $\omega[s_k]$ we
exploit the fact that it is
independent of the smooth part, $v$, of the integrand.
Therefore one may judiciously pick a smooth test function, $g$,  which facilitates the computation of the weight. We choose a test function $g \in C^\infty_c({\mathbb R}^2)$ which is radially symmetric and $g(\zero)=1.$ 
We construct a family of weights $\{\omega_h\}_h$ such that the corrected rule with grid size $h$ integrates exactly our test function $g$: 
\begin{align*}
  & \int_{\mathbb{R}^2}s_k(\mathbf{x})g(\mathbf{x})\text{d}\mathbf{x} = T_{h,\,\mathcal{N}_h}^0[s_k\, g] + h^{k+1}\,\omega_h[s_k]g(\zero)\\  
  \Longrightarrow \hspace{0.4cm} & \omega_h[s_k] := \frac{1}{h^{k+1}}\left[ \int_{\mathbb{R}^2}s_k(\mathbf{x})g(\mathbf{x})\text{d}\mathbf{x}-T_{h,\,\mathcal{N}_h}^0[s_k \,g]\right].
\end{align*}
We define $\omega[s_k]$ by the limit
\begin{equation*}
  \omega[s_{k}] := \lim_{h\to0^+} \omega_h[s_k] = \lim_{h\to0^+} 
  \frac{1}{h^{k+1}}\left[ \int_{\mathbb{R}^2}s_k(\mathbf{x})g(\mathbf{x})\text{d}\mathbf{x}-T_{h,\,\mathcal{N}_h}^0[s_k \,g]\right].
\end{equation*}
Note that since $g$ is chosen to have compact support,
$T_{h,\,\mathcal{N}_h}^0[s_k \,g]$ is a summation of a finite number of terms.
By choosing $g$ radially symmetric, $g(\xx)=g(\vert\xx\vert)$, the two-dimensional Cauchy integral $\int_{{\mathbb{R}}^2}s_k(\xx)g(\xx)\dd\xx$ 
can be efficiently approximated to machine precision,  e.g. using a Gaussian quadrature, by passing to polar coordinates $\xx=r(\cos\theta,\sin\theta)$
\[
\int_{{\mathbb{R}}^2}s_k(\xx)g(\xx)\dd\xx = \int_0^\infty r^{k}g(r)\dd r\,\int_{0}^{2\pi} \phi_k(\cos\theta,\sin\theta)\dd \theta.
\]
%%%%%%%%%%%%%%%%%%%%%%%%%%%%%%%%%%%%%%%%%%%%%%%%

\subsubsection{The case of singular points lying off the grid}\label{sub:first-order-unaligned}
\begin{figure}
    \begin{center}
        \includegraphics[scale=1]{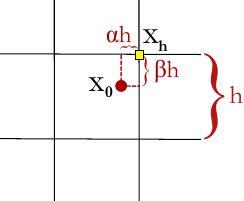} \hspace{2cm}  
        \includegraphics[scale=1]{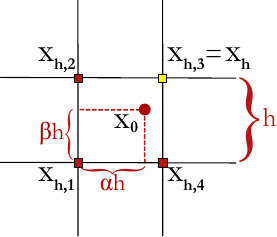}
    \end{center}
    \caption{\textbf{Singularity unaligned to the grid}. The parameters $\alpha,\beta$ are used to characterize the position of the singularity point $\xx_0$ (red circle) relative to the grid in two different settings.
    Left plot (first order correction): position of the singularity point relative to the closest grid node $\xx_h$ (yellow square).
    Right plot (second order correction): position of the singularity point relative to the four surrounding grid nodes $\xx_{h,i}$, $i=1,2,3,4$ (red squares except $\xx_{h,3}=\xx_h$ which is yellow).}
    \label{fig:single_correction_grid_plot}
\end{figure}

In most existing works, see \cite{kapur1997high,aguilar2002high,marin2014corrected,wumartinsson2020zeta}, one assumes that the singularity lies in the origin $\xx_0=\zero$, or equivalently falls in one of the grid nodes.
However, for integrals arising from the IBIM, one must consider the more general case
\begin{equation}\label{eq:sk-integral-x0}
    \int_{{\mathbb{R}}^2} s_k(\xx-\xx_0){v(\xx)}\dd\xx,
\end{equation}
with $s_k$ as in \eqref{eq:sk-integrand} and $\xx_0\notin h\mathbb{Z}^2$. 
We let $\xx_h$ be the grid node closest to $\xx_0$, satisfying
\[
 \xx_{h} = \xx_{h}(\xx_0) = \arg\min_{\xx\in h\mathbb{Z}^2}
 \left\vert\xx-\xx_0\right\vert,\ \ 
 \xx_0=\xx_{h}+(\alpha h, \beta h),\ \
 \alpha,\beta\in\left[-\frac{1}{2},\frac{1}{2}\right),
\]
as shown in the left plot of Figure~\ref{fig:single_correction_grid_plot}.
Correspondingly, we define $\N_{h,1}(\xx_0)=\{\xx_h\}$, and the punctured trapezoidal rule becomes
\[
T_{h,\,\mathcal{N}_{h,1}}^0[f] = h^2\sum_{\xx\in h\mathbb{Z}^2 \setminus \N_{h,1}(\xx_0)}f(\xx).
\]

One can observe from numerical simulations that 
\begin{align}
    \int_{{\mathbb{R}}^2}s_k(\xx-\xx_0)v(\xx)\dd\xx \,=&\, T_{h,\,\mathcal{N}_{h,1}}^0[s_k(\cdot-\xx_0) v(\,\cdot\,)]\nonumber\\
    &\,+h^{k+1}\omega[s_k;\alpha,\beta]v(\xx_h)+\Oo(h^{k+2}). \label{eq:first-corr-term-thm-unaligned}
\end{align}
Hence, the weight $\omega$ 
is still independent of $v$ and $h$, but
now depends
on the relative position of the singularity with respect to the grid, $(\alpha,\beta)$. Moreover,
the function $v$ is evaluated in
$\xx_h$ rather than in the singular point $\xx_0$.
Following this observation, 
we define the \emph{first order correction} $Q_h^1$ to the punctured trapezoidal rule when $\xx_0$ does not fall on the grid as
\begin{equation}
  Q^1_h[s_k(\,\cdot -\xx_0)v(\,\cdot\,)] := 
  T_{h,\,\mathcal{N}_{h,1}}^{0}[s_k(\,\cdot -\xx_0)v(\,\cdot\,)] 
  + h^{k+1}\,\omega[s_k;\alpha,\beta]v(\xx_h).
  \label{eq:single_correction_quadrature}
\end{equation}
Again, this gives an overall error of 
size $\Oo(h^{k+2})$.

The weight only depends on the relative position of the singularity with respect to the grid. 
We therefore set $\xx_0=\zero$ and, fixed $h$ and $(\alpha,\beta)$, shift the grid by $(\alpha,\beta)h$.
 Hence the weight is defined as the limit of the sequence:
\begin{equation}
  \omega[s_k;\alpha,\beta] := \lim_{h\to0^+} \omega_h[s_k;\alpha,\beta]\,,
  \label{eq:single_correction_weight}
\end{equation}
where 
\begin{align*}
  \omega_h[s_k;\alpha,\beta] :=& 
  \dfrac{1}{h^{k+1}}\frac{ \int_{\mathbb{R}^2}{s_k}(\xx){g}(\xx)\text{d}\xx-T_{h,\,\mathcal{N}_{h,1}}^{0}\Big[s_k(\,\cdot -(\alpha,\beta)h) g(\,\cdot -(\alpha,\beta)h)\Big]}{ g( - (\alpha,\beta)h)}.
\end{align*}
The test function $g$ is chosen as in the previous case. The advantages of this choice are going to be the same, e.g. the integral $\int_{\mathbb{R}^2}{s_k}(\xx){g}(\xx)\text{d}\xx=\int_{\mathbb{R}^2}{s_k}(\xx-(\alpha,\beta)h){g}(\xx-(\alpha,\beta)h)\text{d}\xx$ can be computed fast and accurately by passing to polar coordinates.

\subsection{Second order correction}\label{sub:second-order-correction}

The goal now is to build a quadrature rule with error $\Oo(h^{k+3})$ for the integrand $s_k(\xx-\xx_0)v(\xx)$ by approximating
the quadrature error 
$\Rr_h$ of the punctured trapezoidal rule $T_{h,\,\mathcal{N}_{h,1}}^0$
to higher order. 
We have
\begin{equation}\label{eq:error-Rh2}
    \int_{{\mathbb{R}}^2}s_k(\xx-\xx_0)v(\xx)\dd\xx = T_{h,\,\mathcal{N}_{h,1}}^0[s_k(\cdot-\xx_0) v(\,\cdot\,)]
    +\Rr_{h}[s_k(\cdot-\xx_0) v(\,\cdot\,)].
\end{equation}
Since both the singular integral and the punctured trapezoidal rule are linear in $s_k\,v$, we generalize the ansatz for $\Rr_h$
in \eqref{eq:first-corr-term-thm-unaligned} to
achieve higher order accuracy by expanding $v$ {at a} grid node $\tilde\xx_h$ close to $\xx_0$:
\begin{align*}
    \Rr_{h}[s_k(\cdot-\xx_0) v(\,\cdot\,)]
    =& h^{k+1}\omega[s_k;\alpha,\beta]v(\tilde\xx_h) \\
    &+ h^{k+2}
    \left(
    \begin{array}{c}
        u_x[s_k;\alpha;\beta] \\
        u_y[s_k;\alpha;\beta]
    \end{array}
    \right)^T
    \nabla v(\tilde\xx_h)+\Oo(h^{k+3}).
\end{align*}
This ansatz requires three weights $\omega[s_k;\alpha,\beta]\in{\mathbb{R}}$ and $(u_x[s_k;\alpha,\beta],\,u_y[s_k;\alpha,\beta])^T\in\mathbb{R}^2$. 
We further replace the partial derivatives of $v$ at $\tilde\xx_h$ by finite differences of $v$ on
$\{\xx_{h,i}\}_{i=1}^{\tilde p}$, $\xx_{h,1} = \tilde\xx_h$:
\[
\nabla v(\xx_{h,1}) = h^{-1}
\left(
\begin{array}{c}
    \sum_{i=1}^{\tilde p} {\mu_{x,i}}\, v(\xx_{h,i}) \\[0.2cm]
    \sum_{i=1}^{\tilde p} {\mu_{y,i}}\, v(\xx_{h,i})
\end{array}
\right) + \Oo(h),
\]
where $\{\mu_{x,i}\}_{i=1}^{\tilde p}$ are the finite difference weights for the derivative $\frac{\partial}{\partial x} v(\xx_{h,1})$.
Here, the finite differences involve four grid nodes ($\tilde p = 4$) closest to to the singular point $\xx_0$.
They are shown in Figure~\ref{fig:single_correction_grid_plot}
and given by
\begin{equation}
\label{eq:Q4-cutoff}
    \N_{h,4}(\xx_0):=\{\xx_{h,i}\}_{i=1}^4=\{\tilde\xx_h,\ \tilde\xx_h+(0,h),\ \tilde\xx_h+(h,h),\ \tilde\xx_h+(h,0)\}, 
\end{equation}
where $\tilde\xx_h\in h\mathbb{Z}^2$ is the node such that
\[
(\alpha,\beta)=\frac{\xx_0-\tilde\xx_h}{h},~~~\text{for some $\alpha,\beta\in[0,1$).}
\]
We remark that $\alpha,\beta$ are different from the ones for first order correction.

Based on the ansatz above, we define the {\emph{second order correction} $Q_h^2$ to the punctured trapezoidal rule} by: 
\begin{align}\label{eq:Q2-correction}
    Q_h^2[s_k(\,\cdot-\xx_0)\,v(\,\cdot\,)] =& T^0_{h,\,\N_{h,4}}[s_k(\,\cdot-\xx_0)\,v(\,\cdot\,)] \\
    &+ h^{k+1}\sum_{i=1}^4 \omega_i[s_k;\alpha,\beta] v(\xx_{h,i}),\nonumber
\end{align}
where
\[
\omega_i[s_k;\alpha,\beta] :=
\begin{cases}
u_x[s_k;\alpha,\beta]\mu_{x,i} + u_y[s_k;\alpha,\beta]\mu_{y,i} + \omega[s_k;\alpha,\beta], & \text{if } i=1, \\
u_x[s_k;\alpha,\beta]\mu_{x,i} + 
u_y[s_k;\alpha,\beta]\mu_{y,i}, & \text{if } i > 1.
\end{cases}
\]
We note that as long as the finite differences are first order accurate, using them will not change the formal accuracy of the quadrature rule.

The next task is to find a suitable set of weights $\{\omega_i[s_k;\alpha,\beta]\}^4_{i=1}$ for the given $s_k$ and $(\alpha,\beta)$ so that 
\begin{equation}\label{eq:Q2-formula}
    \Rr_h[s_k(\,\cdot-\xx_0)\,v(\,\cdot\,)] = h^{k+1}\sum_{i=1}^4 \omega_i[s_k;\alpha,\beta] v(\xx_{h,i})+\Oo(h^{k+3}).
\end{equation}
As before, the weight only depends on the relative position of the singularity with respect to the grid. We therefore set $\xx_0=\zero$ and, fixed $h$ and $(\alpha,\beta)$, shift the grid by $(\alpha,\beta)h$.
The four closest nodes are then 
\[
\{ \xx_{h,i} \}_{i=1}^4\, = \, \big\{ h(-\alpha,-\beta),\,h(-\alpha,1-\beta),\,h(1-\alpha,1-\beta),\,h(1-\alpha,-\beta) \big\}.
\]
Formula~\eqref{eq:Q2-formula} suggests that we can set up four equations, involving four suitable functions $\{g_j\}_{j=1}^4$: for $j=1,2,3,4$ 
\begin{align*}
\lefteqn{\sum_{i=1}^4 g_j(\xx_{h,i})\,\omega_{i,h}[s_k;\alpha,\beta]
} \hskip 1 cm &
\\
=&\,h^{-k-1}\left(\int_{{\mathbb{R}}^2}s_k(\xx)g_j(\xx)\dd\xx - T^{0}_{h,\,\N_{h,4}}[s_k(\,\cdot-(\alpha,\beta)h)\,g_j(\,\cdot-(\alpha,\beta)h)] \right).
\end{align*}
Fixed $h$ and $(\alpha,\beta)$, this corresponds to imposing that the rule~\eqref{eq:Q2-correction} integrates exactly the functions $s_k(\xx-\xx_0)g_j(\xx-\xx_0)$, $j=1,2,3,4$. Then the weights are found as:
\begin{equation*}
\omega_i[s_k;\alpha,\beta] := \lim_{h\to 0} \omega_{i,h}[s_k;\alpha,\beta] ,\qquad i=1,2,3,4.
\end{equation*}

We choose the test function $g_1=g\in C^\infty_c({\mathbb{R}}^2)$, radially symmetric, such that $g(\zero)=1$ and $\nabla g(\zero)=\zero$. 
This function
behaves like the constant function one {near} $\zero$ and  decays to zero smoothly so that the integrand is compactly supported. These properties facilitate efficient and highly accurate numerical approximation of $\int s_k\,g_j$.
We then use
\begin{align*}
    g_2(x,y) =  x\,g(x,y),\ g_3(x,y) =  y\,g(x,y),\ g_4(x,y) =  x y\,g(x,y),\qquad (x,y)=\xx.
\end{align*}

Out of the four conditions, the first three translate to the weights correctly integrating any function of the type $s_k\, q$, $q\in\mathbb{P}_1$, i.e. $q$ two-dimensional polynomial of degree at most one. 
Three is also the minimum number of points needed in a stencil to have first order accurate $\nabla v$; the fourth node (and consequently the fourth condition) is unnecessary to reach the desired order. It is however useful for it allows us to consider the square four-point stencil \eqref{eq:Q4-cutoff} instead of four different three-point stencils necessary to describe the nodes closest to $\xx_0$.

Computing the right-hand side of the linear system involves evaluating with high accuracy integrals with singular integrands
\[
\int_{{\mathbb{R}}^2}s_k(\xx)g(\xx)\,x^i\,y^j\dd x\dd y = \int_{{\mathbb{R}}^2}\vert\xx\vert^{k-1}\phi_k(\xx/\vert\xx\vert)g(\xx)\,x^i\,y^j\dd x\dd y\ ,\ \ k,i,j\geq0.
\]
By choosing $g(\xx)=g(\vert\xx\vert)$ radially symmetric, we can write the integral in polar coordinates $\xx=r(\cos\theta,\sin\theta)$:
\begin{align*}
&\int_{{\mathbb{R}}^2}\vert\xx\vert^{k-1}\phi_k(\xx/\vert\xx\vert)g(\xx)\,x^i\,y^j\dd x\dd y \\ 
&= \int_{0}^\infty r^{k+i+j}g(r)  \dd r\, \int_{0}^{2\pi} \phi_k(\theta) \cos^i\theta\,\sin^j\theta\,\dd \theta.
\end{align*}
We compute the two factors with high accuracy using Gaussian quadrature. We also reuse the computed values for different parameters $(\alpha,\beta)$.

\subsection{Higher order corrections}
\label{sub:high-order-correction}

We now generalize the approach to construct higher order corrections to the punctured trapezoidal rule for \eqref{thirdordererror}. 
We expand the ansatz \eqref{eq:Q2-formula} used in the previous section to achieve higher order accuracy:
\begin{align*}
\Rr_h[s_k(\,\cdot-\xx_0)\,v(\,\cdot\,)]
=& h^{k+1}\sum_{\vert\nu\vert\leq p-1}h^{\vert\nu\vert}u_{\nu}[s_k;\alpha,\beta]\dfrac{\partial^\nu}{\partial\xx^\nu} v(\xx_{h,1}) + \Oo(h^{k+p+1}),
\end{align*}
where $\nu\in \mathbb{N}^{2}_0$ and the weights $u_\nu\in{\mathbb{R}}$ are independent of $h$ and $v$. This ansatz requires $p(p+1)/2=:p_{\min}$ weights. 
We replace the partial derivatives of $v$ at $\xx_{h,1}$ by sufficiently high order finite differences. Given $\tilde p\geq p_{\min}$, let 
\[
\N_{h,{\tilde p}}(\xx_0):=\{\xx_{h,i}\}_{i=1}^{{\tilde p}}
\]
be a stencil of ${\tilde p}$ nodes close to $\xx_0$, where $\xx_{h,1}$ is such that 
\[
(\alpha,\beta) \, = \, \dfrac{\xx_0-\xx_{h,1}}{h}\,,\ \ \text{for some }\alpha,\beta\in[0,1).
\]
We approximate the derivatives of $v$ using this stencil:
\[
\dfrac{\partial^\nu}{\partial\xx^\nu} v(\xx_{h,1}) = h^{-\vert\nu\vert} \sum_{i=1}^{\tilde p} v(\xx_{h,i})\mu_{\nu,i}+\Oo(h^{ p-\vert\nu\vert}),
\]
where $\{\mu_{\nu,i}\}_{i=1}^{\tilde p}$ are the finite difference weights for the derivative $\frac{\partial^\nu}{\partial \xx^\nu}v(\xx_{h,1})$. 
We finally define the \emph{$p$-th order correction} $Q^p_h$ to the punctured trapezoidal rule as
\begin{equation}\label{eq:Qp-correction}
    Q^p_h[s_k(\cdot\,-\xx_0) v(\,\cdot\,)] := \, T^0_{h,\,\N_{h,{\tilde p}}}[s_k(\,\cdot-\xx_0)v(\,\cdot\,)] + h^{k+1}\sum_{i=1}^{{\tilde p}}\omega_{i}[s_k;\alpha,\beta]v(\xx_{h,i}),
\end{equation}
where
\[
\omega_i[s_k;\alpha,\beta]:=\sum_{\vert\nu\vert\leq p-1} \mu_{\nu,i}\,u_{\nu}[s_k;\alpha,\beta] , \qquad i=1,\dots,\tilde p\,.
\]
As long as the finite differences for $\frac{\partial^\nu}{\partial \xx^\nu}v(\xx_{h,1})$ have error $\sim \Oo(h^{p-\vert\nu\vert})$ they will not affect the formal accuracy of the quadrature rule.

We now have to find a suitable set of weights $\{\omega_i[s_k;\alpha,\beta]\}^{\tilde p}_{i=1}$ for the given $s_k$ and $(\alpha,\beta)$ so that for any smooth function $v$
\begin{equation}\label{eq:Qp-formula}
    \Rr_h[s_k(\,\cdot-\xx_0)\,v(\,\cdot\,)] = h^{k+1}\sum_{i=1}^{\tilde p} \omega_i[s_k;\alpha,\beta] v(\xx_{h,i})+\Oo(h^{k+1+p}).
\end{equation}

Analogously to Section~\ref{sub:second-order-correction}, the weights only depend on the relative position of the singularity with respect to the grid. We therefore set $\xx_0=\zero$ and, fixed $h$ and $(\alpha,\beta)$, shift the grid by $(\alpha,\beta)h$.\\
Formula~\eqref{eq:Qp-formula} suggests that we may set up $\tilde p$ equations, involving $\tilde p$ suitable test functions $\{g_j\}_{j=1}^{\tilde p}$, to uniquely define the weights $\{\omega_{i}\}_{i=1}^{\tilde p}$. We proceed as in the previous Section and define the family of weights $\{\omega_{i,h}\}_{i=1}^{\tilde p}$ solution to
\begin{align}
\lefteqn{\sum_{i=1}^{\tilde p} g_j(\xx_{h,i})\,\omega_{i,h}[s_k;\alpha,\beta]}
\hskip 1 cm &\label{eq:def-w_ih}\\
=\,& h^{-k-1}\left(\int_{{\mathbb{R}}^2}s_k(\xx)g_j(\xx)\dd\xx - T^{0}_{h,\,\N_{h,\tilde p}}[s_k(\,\cdot-(\alpha,\beta)h)\,g_j(\,\cdot-(\alpha,\beta)h] \right),\nonumber
\end{align}
for $j=1,\dots,\tilde p$.
Fixed $h$ and $(\alpha,\beta)$, this corresponds to imposing that the rule~\eqref{eq:Qp-correction} integrates exactly the functions $s_k(\xx-\xx_0)g_j(\xx-\xx_0)$, $j=1,\dots,\tilde p$. 
Then the weights are found as:
\begin{equation}\label{eq:weights-Qp-limit}
\omega_i[s_k;\alpha,\beta] := \lim_{h\to 0} \omega_{i,h}[s_k;\alpha,\beta]\ ,\ \ i=1,\dots,\tilde p.
\end{equation}

We use the function $g$ similar to the one considered in Section~\ref{sub:second-order-correction}, with the additional conditions that $\frac{\partial^\nu}{\partial\xx^\nu}g(\zero)=0$ for all $\vert\nu\vert\leq p-1$. 
This ensures that $g$ is similar enough to the constant function $g\equiv 1$ near $\zero$.\\

By choosing the $p_{\min}$ functions $\{ g_j\}_{j=1}^{p_{\min}}$ equal to $g$ multiplied by the $p_{\min}$ monomials of degree at most $p-1$ ($x^i y^j$, $i,j\in\mathbb{N}_0$, $i+j\leq p-1$) we impose that the method \eqref{eq:Qp-correction} integrates exactly all integrands of the type $s_k\, q$, $q\in\mathbb{P}_{p-1}$, i.e. two-dimensional polynomials of degree at most $p-1$. 
The additional $\tilde p-p_{\min}$ functions can be chosen for example as $g$ multiplied by two-dimensional monomials of degree higher than $p$.\\

We use $\tilde p\geq p_{\min}$ because $p_{\min}$ may not fit well with standard stencils. Thus it is possible to use more nodes than $p_{\min}$ and impose additional conditions. For example in our implementations for first, second, third, and fourth order corrections we used $\tilde p\geq p_{\min}$ as shown in Table~\ref{tab:p-bar-stencils}. A visualization of these stencils can be seen in Figure~\ref{fig:all_correction_grid_plot}.

\begin{table}
    \caption{\textbf{Correction order and corresponding correction nodes}. To increase the order of accuracy by $p$, the minimum number of nodes to correct is $p_{\min}=p(p+1)/2$ but more nodes can be used. $\tilde p\geq p_{\min}$ is the number of nodes we used in our tests, corresponding to the stencils showed in Figure~\ref{fig:all_correction_grid_plot}.}
    \label{tab:p-bar-stencils}
    \begin{center}
    \begin{tabular}{c|c|c}
        \ $p$\ & \ $p_{\min}$\ & \ $\tilde p$\ \\[0.1cm]
        \hline
        1 & 1 & 1 \\ 
        2 & 3 & 4 \\
        3 & 6 & 6 \\
        4 & 10 & 12
    \end{tabular}
    \end{center}
\end{table}

\subsection{High order quadratures}
\label{sub:additive-splitting}

In the previous sections
we have shown 
how to deal with integrands of the kind \eqref{eq:sk-integrand}
\[
s_k(\xx-\xx_0)v(\xx)=\vert\xx-\xx_0\vert^{k-1}\phi_k((\xx-\xx_0/\vert\xx-\xx_0\vert)v(\xx) ,\qquad k=0,1,2,\dots\,,
\]
wherever the singularity point $\xx_0$ may lie, which means we can correct the trapezoidal rule for all terms in the expansion \eqref{eq:s-expansion} 
\begin{align*}
    s(\xx-\xx_0) =&\, \dfrac{1}{\vert\xx-\xx_0\vert}\phi_0\left(\frac{\xx-\xx_0}{\vert\xx-\xx_0\vert}\right) + \phi_1\left(\frac{\xx-\xx_0}{\vert\xx-\xx_0\vert}\right) \\ &+ \vert\xx-\xx_0\vert\,\phi_2\left(\frac{\xx-\xx_0}{\vert\xx-\xx_0\vert}\right) + \cdots
\end{align*}
of the singular function \eqref{eq:s-singular-integrand}.
If we know these terms explicitly we can build a high order corrected trapezoidal rule for the integral
\[
\int_{{\mathbb{R}}^2} s(\xx-\xx_0)v(\xx)\dd\xx.
\]
We demonstrate the idea of successive corrections
by deriving a second and then a third order accurate quadrature rule. 

We first write
\begin{align*}
    s(\xx-\xx_0)v(\xx) =& s_0(\xx-\xx_0)v(\xx)+\left( s(\xx-\xx_0)- s_0(\xx-\xx_0) \right){v(\xx)} \\
    =& s_0(\xx-\xx_0)v(\xx)+\triangle_0 s(\xx-\xx_0)v(\xx).
\end{align*}
Lemma \ref{lem:Rq-remainder-trapezoidal-rule} states that the punctured trapezoidal rule is second order accurate for integrating $\triangle_0 s$. If we apply the first order correction \eqref{eq:single_correction_quadrature} to the punctured trapezoidal rule for the first term, we get a second order approximation.
The explicit formula, with 
$\N_{h,1}=\{\xx_h\}$ as in Section \ref{sub:first-order-unaligned} and relative grid shifts $(\alpha_1,\beta_1)$, is:
\begin{align*}
    \mathcal{U}^{2}_{h}[s(\,\cdot-\xx_0)v(\,\cdot\,)] :=& \, Q_h^1\left[s_0(\,\cdot-\xx_0)v(\,\cdot\,) \right] + T^0_{h,\,\N_{h,1}}[\Delta_0(\,\cdot-\xx_0)v(\,\cdot\,)] \\[0.2cm]
    =&\,  h^2\sum_{\xx\in h\mathbb{Z}^2 \setminus \N_{h,1}(\xx_0)}s(\xx-\xx_0)v(\xx) + h\,\omega[s_0;\alpha_1,\beta_1]\,v(\xx_h).
\end{align*}
This was the approach used in \cite{izzo2021corrected}, although there \eqref{eq:single_correction_quadrature} was used also on the second term instead of the punctured trapezoidal rule.

To achieve third order, we expand $s$ further:
\begin{align}
    s(\xx-\xx_0)v(\xx) =&\, s_0(\xx-\xx_0)v(\xx) + s_1(\xx-\xx_0)v(\xx)  \nonumber\\[0.1cm]
    &+ \left[s(\xx-\xx_0)-s_0(\xx-\xx_0)-s_1(\xx-\xx_0)\right]v(\xx)\nonumber\\[0.1cm]
    =&\, s_0(\xx-\xx_0)v(\xx) + s_1(\xx-\xx_0)v(\xx)+\triangle_1 s(\xx-\xx_0)v(\xx)\label{eq:additive-splitting}.
\end{align}
We then use the second order correction \eqref{eq:Q2-correction} for integrating the first term, first order correction \eqref{eq:single_correction_quadrature} for integrating the second, and the (uncorrected) punctured trapezoidal rule for $\triangle_1 s$;
by Lemma \ref{lem:Rq-remainder-trapezoidal-rule} it is third order accurate for $\triangle_1s$.

We use the set of correction nodes $ \N_{h,4}(\xx_0) = \{\xx_{h,i}\}_{i=1}^4$ and define the corresponding relative grid shift $(\alpha_2,\beta_2)$.
Then the third order accurate rule $\mathcal{U}^3_h$ is
\begin{align}
    \mathcal{U}^{3}_{h}[s(\cdot-\xx_0)v(\,\cdot\,)] :=& \,  Q_h^2\left[s_0(\,\cdot-\xx_0){v(\,\cdot\,)} \right] +  Q_h^1\left[ s_1(\,\cdot-\xx_0){v(\,\cdot\,)} \right] \nonumber\\[0.2cm]
    &+ T^0_{h,\,\N_{h,1}}[\triangle_1 s(\,\cdot\,-\xx_0){v(\,\cdot\,)}] \label{eq:Q-gen-third-order-2D}\\[0.15cm]
    =&  \, h^2\sum_{\xx\in h\mathbb{Z}^2 \setminus \N_{h,4}(\xx_0)}s(\xx-\xx_0)v(\xx) + h\sum_{i=1}^4 \omega_i[{s_0};\alpha_2,\beta_2] v(\xx_{h,i}) \nonumber\\[0.2cm] 
    &+ h^2\,\omega[s_1;\alpha_1,\beta_1]\, v(\xx_{h})  \nonumber\\[0.2cm]
    & + h^2\sum_{\xx\in \N_{h,4}(\xx_0)\setminus \N_{h,1}(\xx_0)}  \big\{ s(\xx-\xx_0)-s_0(\xx-\xx_0) \big\}v(\xx).\nonumber
\end{align}

\begin{figure}
    \begin{center}
        \includegraphics[scale=1]{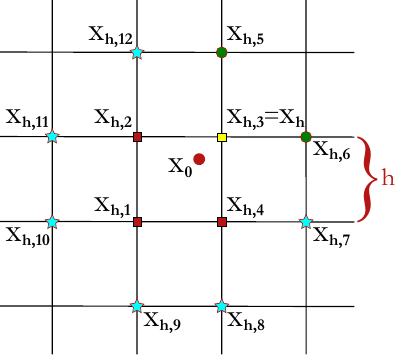}
    \end{center}
    \caption{\textbf{Example of correction stencils}. The stencils we tested for corrections $p=1$ ($\tilde p=1$ node: yellow square), $p=2$ ($\tilde p=4$ nodes: yellow and red squares), $p=3$ ($\tilde p=6$ nodes: yellow and red squares, and green circles), and $p=4$ ($\tilde p=12$ nodes: yellow and red squares, green circles, and cyan stars). The singularity node is $\xx_0$ (red circle). The nodes are $\{\xx_{h,i}\}_{i=1}^{12}$, and $\N_{h,\tilde p}=\{\xx_{h,i}\}_{i=1}^{\tilde p}$ except for $\N_{h,1}=\{ \xx_h\}$. }
    \label{fig:all_correction_grid_plot}
\end{figure}

In general, given the singular function $s(\xx-\xx_0)v(\xx)$, in order to build a quadrature rule $\mathcal{U}^p_h$ of order $p\ge 2$ we need explicitly the first $p-1$ ($k=0,\dots,p-2$) terms of the expansion \eqref{eq:s-expansion}
\[
s(\xx) = \sum_{k=0}^{p-2} s_k( {\xx}) + \triangle_{p-2}s(\xx) = \sum_{k=0}^{p-2} \vert\xx\vert^{k-1}\phi_k\left( \dfrac{\xx}{\vert\xx\vert} \right) + \triangle_{p-2}s(\xx),
\]
and apply to the term $s_k(\xx)$ the $(p-k-1)$-th order correction $Q^{p-k-1}_h$ to the trapezoidal rule. The punctured trapezoidal rule is used for $\triangle_{p-2}s$. 
\begin{equation}\label{eq:Qp-general-form}
    \mathcal{U}^p_h[s(\,\cdot-\xx_0) v(\,\cdot\,)] := \sum_{k=0}^{p-2} Q_h^{p-1-k}\left[ s_k({\,\cdot-\xx_0})v(\,\cdot\,) \right] + T_{h,\,\N_{h,1}}^0[\triangle_{p-2}s(\,\cdot-\xx_0)v(\,\cdot\,)].
\end{equation}
We can find an explicit expression for the quadrature rule $\mathcal{U}_h^p$ by specifying the stencils we use for the correction nodes. We denote by $\N_{n,\tilde p(p)}$ the stencil of $\tilde p(p)$ correction nodes to increase the order by $p$. We assume that the stencils are increasing: $\N_{n,\tilde p(p)}\subset \N_{n,\tilde p(p+1)}$.
For example in our tests we took $\tilde p(1)=1$, $\tilde p(2)=4$, $\tilde p(3)=6$, $\tilde p(4)=12$, and 
$\N_{h,1}(\xx_0):=\{\xx_h\}=\{\xx_{h,3}\}$, $\N_{h,4}(\xx_0)=\{\xx_{h,i}\}_{i=1}^4$, $\N_{h,6}(\xx_0)=\{\xx_{h,i}\}_{i=1}^6$, $\N_{h,12}(\xx_0)=\{\xx_{h,i}\}_{i=1}^{12}$, so that $\N_{h,1}\subset \N_{h,4}\subset \N_{h,6}\subset \N_{h,12}$. 
This is shown in Table~\ref{tab:p-bar-stencils} and Figure~\ref{fig:all_correction_grid_plot}.
We call $\alpha_p,\beta_p$ the parameters describing the shift with respect to $\xx_0$ of the stencil of $\tilde p(p)$ nodes:
\begin{align}
    \mathcal{U}^p_h[s(\,\cdot-\xx_0) v(\,\cdot\,)] =&\, h^2\sum_{\xx\in h \mathbb{Z}^2 \setminus \N_{h,\tilde p(p-1)}} s(\xx-\xx_0)v(\xx) \nonumber\\[0.2cm]
    & + h^2 \sum_{k=1}^{p-3}\, \sum_{\xx\in \N_{h,\tilde p(p-1)}\setminus \N_{\tilde p(p-k-1)} }   s_k({\xx-\xx_0})v(\xx) \label{eq:Qp-general-explicit}\\[0.2cm]
    & + h^{p-1}\,\omega [ s_{p-2}; \alpha_{1},\beta_{1} ]v(\xx_{h}) \nonumber\\[0.2cm]
    & + \sum_{k=0}^{p-3} h^{k+1} \sum_{i=1}^{\tilde p(p-k-1)} \omega_i [ s_k; \alpha_{p-k-1},\beta_{p-k-1} ] v(\xx_{h,i}) \nonumber\\[0.2cm]
    & + h^2 \sum_{\xx\in \N_{h,\tilde p(p-1)}\setminus \N_{h,\tilde p(1)}} \left\{ s(\xx-\xx_0)-\sum_{k=0}^{p-3}s_k({\xx-\xx_0}) \right\}v(\xx) \nonumber .
\end{align}

In Section \ref{sec:num_2D} we show tests for the quadrature method \eqref{eq:Qp-general-explicit} by combining first, second, third, and fourth order corrections.

\subsection{Approximation and tabulation of the weights}
Given functions of the kind \eqref{eq:sk-integrand},
$s_k(\xx)=\vert\xx\vert^{k-1}\phi_k(\xx/\vert\xx\vert)$, we want to compute the weights $\{\omega_i[s_k;\alpha,\beta]\}_{i=1}^{\tilde p}$, defined by
\eqref{eq:def-w_ih} and
\eqref{eq:weights-Qp-limit}. 

Fixed $k\ge 0$ and $(\alpha,\beta)$, we write the function $s_k$ (specifically its factor $\phi_k$) using its Fourier series:
\begin{align*}
\xx\,=&\,\vert\xx\vert\big(\cos(\psi(\xx)),\sin(\psi(\xx))\big), \\[0.3cm]
s_k(\xx)\,=&\,\vert\xx\vert^{k-1}\phi_k(\xx/\vert\xx\vert)=\vert\xx\vert^{k-1}\phi_k(\psi(\xx))\\
=&\,\vert\xx\vert^{k-1}\left(a_0+\sum_{j=1}^\infty \big( a_j\cos(j\psi(\xx))+b_j\sin(j\psi(\xx)) \big)\right),
\end{align*}
where $\{a_j\}_{j=0}^\infty$ and $\{b_j\}_{j=1}^\infty$ are the Fourier coefficients
of $\phi_k$.
Then, by linearity of the weights with respect to $s_k$, we can write them as
\begin{align*}
\omega_i[s_k;\alpha,\beta] \,=&\, a_0\,\omega_i\left[\vert\xx\vert^{k-1};\alpha,\beta\right]\\
&\, + \sum_{j=1}^\infty \Big(a_j\,\omega_i\left[\vert\xx\vert^{k-1}\cos(j\psi(\xx));\alpha,\beta\right]\\
&\,\qquad +b_j\,\omega_i\left[\vert\xx\vert^{k-1}\sin(j\psi(\xx));\alpha,\beta\right]\Big),
\end{align*}
with $i=1,\dots,\tilde p$. 
We can then approximate and tabulate the weights 
$\omega_i[s_k;\alpha,\beta]$
in the following way. We fix a stencil of parameters 
$\{(\alpha_m,\beta_n)\}_{m,n}$
around $(\alpha,\beta)$ 
and basis functions $\{c_{m,n}(\alpha,\beta)\}_{m,n}$
such that we can approximate a function $f:{\mathbb{R}}^2\to{\mathbb{R}}$ in $(\alpha,\beta)$ as
\[
f(\alpha,\beta)\approx \sum_{m,n}c_{m,n}(\alpha,\beta)\,f(\alpha_m,\beta_n).
\]
We let $N$ be
the number of Fourier modes used to approximate the weights.  
Then, given $\phi_k$, we first find the $2N+1$ coefficients 
$a_0,\{a_j,b_j\}_{j=1}^N$ 
by using the Fast Fourier Transform.
Then, 
\begin{align*}
\omega_i[s_k;\alpha_m,\beta_n] \,\approx\,& a_0\,\omega_i\left[\vert\xx\vert^{k-1};\alpha_m,\beta_n\right] \\ 
& + \sum_{j=1}^N \Big( a_j\,\omega_i\left[\vert\xx\vert^{k-1}\cos(j\psi(\xx));\alpha_m,\beta_n\right] \\ &\qquad + b_j\,\omega_i\left[\vert\xx\vert^{k-1}\sin(j\psi(\xx));\alpha_m,\beta_n\right]\Big),
\end{align*}
and we can approximate the weight for $(\alpha,\beta)$ via
\[
\omega_i[s_k;\alpha,\beta] \approx \sum_{m,n}c_{m,n}(\alpha,\beta)\,\omega_i[s_k;\alpha_m,\beta_n],\qquad i=1,\dots,\tilde p.
\]
So, for all expansion terms $k=0,1,\dots,p-2$ used in \eqref{eq:Qp-general-form}, and the corresponding corrections  $Q^{p-1-k}_h$, we need to compute and store
the weights for the 
following constant and trigonometric functions,
\[
\left.
\begin{array}{l}
\omega_i\left[\vert\xx\vert^{k-1};\alpha_m,\beta_n\right]\\[0.15cm]
\omega_i\left[\vert\xx\vert^{k-1}\cos(j\psi(\xx));\alpha_m,\beta_n\right] \\[0.15cm] 
\omega_i\left[\vert\xx\vert^{k-1}\sin(j\psi(\xx));\alpha_m,\beta_n\right]%\ \ j=1,\dots,N,\ \ i=1,\dots,\tilde p(p-1-k),\ \ m,n\in\mathcal{I}.
\end{array}
\right\} \
\begin{array}{l}
    j=1,\dots,N,\\
    i=1,\dots,\tilde p(p-1-k),\\
    \text{and all $m,n$ in the stencil for $(\alpha,\beta)$}.
\end{array}
\]
\begin{remark}
The weights $\omega_i[s_k;\alpha_m,\beta_n]$ are formally the limits of $\omega_{i,h}[s_k;\alpha_m,\beta_n]$ defined in \eqref{eq:def-w_ih}. We approximate the limit by $\omega_{i,h^*}$, where
\[
h^*:=2^{-M}, \ \ M := \arg \min_{j=1,2,3,\dots} \left\{ \vert\omega_{i,2^{-j}}-\omega_{i,2^{-j-1}}\vert\leq \text{Tol}\right\}.
\]
In the simulations presented in Section~\ref{sec:numerical-tests}, to compute $\{\omega_{i}\}_{i=1}^{\tilde p(p)}$, we use $\text{Tol}=10^{-8}$ for $p=1,2,3$, and $\text{Tol}=10^{-4}$ for $p=4$.
\end{remark}

\section{Evaluating layer potentials in the implicit boundary integral formulation}
\label{sec:application-layer-potentials-3D}

We apply the 
high order quadrature methods from
Section~\ref{sec:corrected-trapezoidal-rules} to layer potentials used
in Implicit Boundary Integral Methods (IBIM).
To make the exposition clear
we adopt the following convention.
\begin{notation}
We distinguish between variables in ${\mathbb{R}}^2$ and ${\mathbb{R}}^3$ by using boldface variables for vectors in ${\mathbb{R}}^2$ and boldface variables with a bar for vectors in ${\mathbb{R}}^3$. For example, $\xx\in{\mathbb{R}}^2$ and $\bar\xx\in{\mathbb{R}}^3$.\\
Moreover, for a vector ${\mathbf{y}}=(y_1,y_2)\in {\mathbb{R}}^2$ 
and scalar $y_3\in{\mathbb{R}}$ we frequently write $\bar{\mathbf{y}}=({\mathbf{y}},y_3)$ to mean the vector $(y_1,y_2,y_3)\in{\mathbb{R}}^3$.
For example, 
when
$f:{\mathbb{R}}^3\to{\mathbb{R}}^m$, we use the 
notations
$f(\bar{\mathbf{y}})\equiv f({\mathbf{y}},y_3)\equiv f(y_1,y_2,y_3)$.
\end{notation}
We consider the general form of a layer potential on a smooth, closed and bounded surface $\Gamma\subset {\mathbb R}^3$,
\begin{equation}\label{eq:general_S}
\int_{\Gamma}K(\bar\xx^*,\bar{\mathbf{y}})\rho(\bar{\mathbf{y}})\dd\sigma_{\bar{\mathbf{y}}} , \qquad \bar\xx^*\in\Gamma,
\end{equation}
with $K$  defined by one of the following kernels:
\begin{equation}
\begin{array}{rrrl}
\text{(single-layer, SL)}:& G_{0}(\bar\xx^*,\bar{\mathbf{y}}) = & \dfrac{1}{4\pi} & \dfrac{1}{\vert\bar\xx^*-\bar{\mathbf{y}}\vert},\\ \text{(double-layer, DL)}:& \dfrac{\partial G_{0}}{\partial \nn_{y}}(\bar\xx^*,\bar{\mathbf{y}}) = & \dfrac{1}{4\pi} & \dfrac{(\bar\xx^*-\bar{\mathbf{y}})^T \nn_{y}}{\vert\bar\xx^*-\bar{\mathbf{y}}\vert^{3}},\\
\text{(double-layer conjugate, DLC)}:& \dfrac{\partial G_{0}}{\partial \nn_{x}}(\bar\xx^*,\bar{\mathbf{y}}) = & -\dfrac{1}{4\pi} & \dfrac{(\bar\xx^*-\bar{\mathbf{y}})^T \nn_{x}}{\vert\bar\xx^*-\bar{\mathbf{y}}\vert^{3}}.
\end{array}\label{eq:layer-kernels}
\end{equation}
In \eqref{eq:layer-kernels}, the vector $\nn_{x}$ is the normal vector to $\Gamma$ at $\bar\xx^*$, pointing into the unbounded region ${\mathbb{R}}^3\setminus \overline{\Omega}$, where $\Omega$ is the bounded region enclosed by $\Gamma$. Analogously $\nn_y$ is the normal vector to $\Gamma$ at $\bar {\mathbf{y}}$. 
In preparation for the
formulation of 
the implicit boundary
integral methods we first define
$d_\Gamma:{\mathbb R}^3\mapsto {\mathbb R}$ 
to be the signed distance to the surface such that $d_\Gamma$ is negative inside $\Omega$. Moreover,
we let $P_\Gamma:{\mathbb{R}}^3\to\Gamma$ be the closest point
mapping that takes $\bar{\mathbf{y}}$ to a closest point on $\Gamma$:
\[
    P_\Gamma(\bar{\mathbf{y}}) \in \arg\min_{\bar{\mathbf{z}}\in\Gamma} \vert\bar{\mathbf{y}}-\bar{\mathbf{z}}\vert^2.
\]
Let $\mathcal{C}_\Gamma\subset\mathbb{R}^3$ be the set containing the all the points that have non-unique closest points on $\Gamma$.
The \emph{reach} $\tau_{\Gamma}$ of $\Gamma$ is 
defined as 
\[
\tau_{\Gamma} := \inf_{\bar{\mathbf{x}}\in \Gamma,\, \bar{\mathbf{y}}\in\mathcal{C}_\Gamma} \vert\bar{\mathbf{x}}-\bar{\mathbf{y}}\vert.
\]
It depends on the local geometry (the curvatures) and the global structure of $\Gamma$ (the Euclidean and geodesic distances between any two points on $\Gamma$).
The reach is positive $\tau_\Gamma>0$ if $\Gamma$ is $C^{1,\alpha}$ for some $\alpha>0$. 
The closest point mapping $P_\Gamma$ is invertible in the tubular neighborhood 
\[
T_\varepsilon:=\{\bar\xx\in{\mathbb{R}}^3\,:\,\vert d_\Gamma(\bar\xx)\vert\leq \varepsilon\}\subset {\mathbb{R}}^3,~~~\varepsilon<\tau_\Gamma.
\]
In this paper, we will assume that $\Gamma$ is a closed bounded $C^2$ surface so that the mean and Gaussian curvatures are defined everywhere on the surface.  
Consequently, when $\bar{\mathbf{y}}$ lies within the reach of $\Gamma$, 
we have the explicit formula
\[
P_\Gamma(\bar{\mathbf{y}}) = \bar{\mathbf{y}}-d_\Gamma(\bar{\mathbf{y}})\nabla d_\Gamma(\bar{\mathbf{y}}).
\]
The surface integral (\ref{eq:general_S}) 
can then be reformulated into an equivalent volume integral using the Implicit Boundary Integral Methods  \cite{kublik2013implicit,tsaikublik16},
\begin{equation}\label{eq:volume_potential}
I_\varepsilon [\rho](\bar\xx)= \int_{T_\varepsilon} {K}(\bar\xx,P_\Gamma(\bar{\mathbf{y}}))\rho(P_\Gamma(\bar{\mathbf{y}})) \delta_{\Gamma, \varepsilon}(\bar{\mathbf{y}}) \dd\bar{\mathbf{y}}.  
\end{equation}
The ``delta'' function is defined as $\delta_{\Gamma,\varepsilon}(\bar{\mathbf{y}}):=\delta_\varepsilon(d_\Gamma(\bar{\mathbf{y}}))J_{d_\Gamma(\bar{\mathbf{y}})}(\bar{\mathbf{y}})$, where  
\[
J_\eta(\bar{\mathbf{y}})=1+2\eta H(\bar{\mathbf{y}})+\eta^2 G(\bar{\mathbf{y}}),
\]
with $H(\bar{\mathbf{y}})$ and $G(\bar{\mathbf{y}})$ denoting respectively the mean and Gaussian curvatures of $\Gamma_\eta:=\{\bar\xx\in{\mathbb{R}}^3\,:\,d_\Gamma(\bar\xx)=\eta\}$ (see e.g. \cite{tsaikublik16}). For $\vert\eta\vert<\tau_\Gamma$, $J_\eta$ is bounded away from zero. Moreover,  
\[
\delta_\varepsilon(\eta)=\frac{1}{\varepsilon}\delta\left(\frac{\eta}{\varepsilon}\right)
\]
is smooth compactly supported in $(-\varepsilon,\varepsilon)$ with unit mass.
This is achieved by using $\delta\in C^\infty_c({\mathbb{R}})$ compactly supported in $(-1,1)$ with $\int_{\mathbb{R}} \delta(\eta)\dd\eta=1$.

It turns out that for any positive $\varepsilon$, smaller than the reach of $\Gamma$,
the IBIM is equal to the original
layer potential for all $\bar\xx\in{\mathbb R}^3$,
\begin{equation}\label{eq:equivalence}
I_\varepsilon[\rho](\bar\xx) \equiv \int_{\Gamma}K(\bar\xx,\bar{\mathbf{y}})\rho(\bar{\mathbf{y}})\dd\sigma_{{\mathbf{y}}},~~~\bar\xx\in\mathbb{R}^3.
\end{equation}
If the surface $\Gamma$ is smooth, the closest point  mapping is also smooth (\cite{delfour2011shapes}, Ch.~7, \textsection 8, Thm 8.4). As a consequence, if $\rho$ is a smooth function on $\Gamma$ and $\Gamma$ is smooth, then 
{$\rho(P_\Gamma(\bar{\mathbf{y}}))\delta_{\Gamma, \varepsilon}(\bar{\mathbf{y}})$ is a smooth
function ${\mathbb R}^3$,
compactly supported in $T_\varepsilon$.}

\subsection{Singular integrand in three dimensions and correction plane by plane}
\label{sec:plane-by-plane}

In this section, we will construct high order quadratures 
for $I_\varepsilon[\rho](\bar\xx^*)$
in \eqref{eq:volume_potential}
from the two-dimensional corrected trapezoidal rules \eqref{eq:Q-gen-third-order-2D}.
The three dimensional quadrature rules will be defined
as the sum of 
integration over different coordinate planes, 
where
the two dimensional corrected trapezoidal rule 
from the previous Section
is applied to approximate the integration over each plane. The particular selection of the coordinate planes depends on the normal vector of $\Gamma$. 

\begin{figure}
    \begin{center}
        \includegraphics[height=3cm]{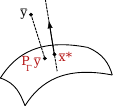}\hspace{0.2cm}\includegraphics[scale=1]{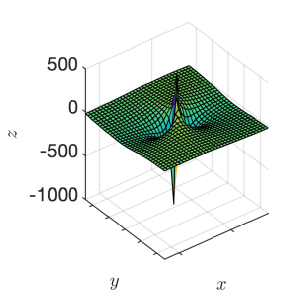}\includegraphics[scale=1]{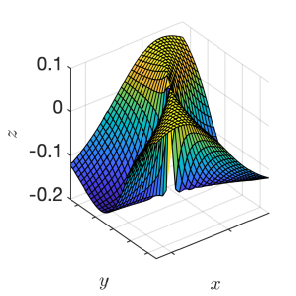}
    \end{center}
    \caption{\textbf{IBIM kernel singular behavior}. The kernel $K(\bar\xx^*,P_\Gamma(\bar{\mathbf{y}}))$, for fixed $\bar\xx^*\in\Gamma$ and $\bar{\mathbf{y}}\in T_\varepsilon$, is singular along the normal $\nn$ to $\bar\xx^*$. The left figure illustrates how the kernel becomes singular for $\bar{\mathbf{y}}$ approaching any point of the line passing through $\bar\xx^*$ with direction $\nn$.
    The center plot shows the double-layer conjugate kernel $K=\frac{\partial G}{\partial \n_x}$ 
    plotted on a plane with fixed $z$, $\bar{\mathbf{y}}\in \{(x,y,z)\,:\,x,y\in{\mathbb{R}}\}$. The function will then have a point singularity in $\yy_0(z)$, and we plot the kernel for $(x,y)$ close to $\yy_0(z)$.
    The right most plot shows the same kernel, multiplied by $\vert{\mathbf{y}}-\yy_0(z)\vert$.
    }
    \label{fig:singularity_along_normals}
\end{figure}

Without loss of generality, we consider
a target point,
$\bar\xx^*=(x^*,y^*,z^*)\in\Gamma$, at which
the surface normal is $\nn=(n_1,n_2, 1)$.
The way to treat other cases are explained in Section~\ref{sub:change-coordinates}.
We denote the integrand 
in \eqref{eq:volume_potential}
by $f$,
\begin{equation}\label{eq:kernelfunction}
f(\bar{\mathbf{y}}):={K}(\bar\xx^*, P_\Gamma(\bar{\mathbf{y}}))\rho(P_\Gamma(\bar{\mathbf{y}}))\delta_{\Gamma,\varepsilon}(\bar{\mathbf{y}}),\qquad
\bar{\mathbf{y}}=(x,y,z).
\end{equation}
To approximate \eqref{eq:volume_potential}
the standard trapezoidal rule is 
first applied in the $z$-direction.
With the grid points
$z_k=kh$, we get
\begin{equation}\label{eq:plane-by-plane-sum}
I_\varepsilon[\rho](\bar\xx^*)=
\int_{\mathbb{R}^{3}}f(x,y,z)\dd x\dd y\dd z 
\approx 
h\sum_{k}\int_{\mathbb{R}^{2}}f(x,y,z_k)\dd x\dd y,
\end{equation}
where we used the
fact that $f$ is compactly supported
in $T_\varepsilon$.
We note that 
$f$ is singular along the line
\begin{equation}\label{eq:y0(z)_expression}
\bar\yy_0(z)=\bar\xx^*+(z-{z}^*)\nn,
\qquad z\in{\mathbb R},
\end{equation}
since $P_\Gamma(\bar\yy_0(z))=\bar\xx^*$
for all $z$.
Therefore,
$f(\cdot,\cdot,z)$ is singular at one point
for each fixed $z$, by the assumption on $\nn$.
Below we will derive the
form of this singularity and we will
show that it is of the same type
\eqref{eq:s-singular-integrand}
as considered in Section~\ref{sec:corrected-trapezoidal-rules}. See Figure~\ref{fig:singularity_along_normals} for an illustration of the line singularity, and an example of the singular behavior.
Hence,
we can use the corrected trapezoidal rules
to approximate each
integral in the sum in
\eqref{eq:plane-by-plane-sum}.

To connect back to the notation
in Section~\ref{sec:corrected-trapezoidal-rules} we 
write $\bar{\mathbf{y}}=({\mathbf{y}},z)$ and
$\bar\yy_0(z)= (\yy_0(z),z)$.
Then we
factorize $f$, for a fixed $z$, as
\begin{equation}\label{eq:f_s_v_3D_splitting}
f({\mathbf{y}}, z) = 
s({\mathbf{y}}-\yy_0(z);z) v({\mathbf{y}}, z).
\end{equation}
where
\begin{equation}\label{eq:svdefine}
s({\mathbf{y}};z)={K}(\bar\xx^*, P_\Gamma({\mathbf{y}}+\yy_0(z), z)),
\qquad
v({\mathbf{y}},z) =\rho(P_\Gamma({\mathbf{y}}, z))\delta_{\Gamma,\varepsilon}({\mathbf{y}}, z).
\end{equation}
Note that the type of singularity for $s$
depends on the properties of $\Gamma$ at the target point (such as principal curvatures, principal directions, normal).
Moreover, $s$ depends smoothly on $z$. 

We then use the corrected trapezoidal rule $\mathcal{U}^3_h[f]$ \eqref{eq:Q-gen-third-order-2D} to compute the integrals on each plane,
\[
\int_{\mathbb{R}^2}f({\mathbf{y}}, z_k)\dd{\mathbf{y}} =
\int_{\mathbb{R}^2}
s({\mathbf{y}}-\yy_0(z_k);z_k) v({\mathbf{y}}, z_k) \dd {\mathbf{y}}
\,\approx\, \mathcal{U}^3_h[
s(\,\cdot-\yy_0(z_k);z_k) v(\,\cdot\,, z_k)].
\]
We denote by ${{\mathbf{y}}}_{\Delta}(z)$ and $(\alpha_1(z),\beta_1(z))$ 
the closest grid node to $\yy_0(z)$ and the relative grid shift parameters
respectively, as defined in Section~\ref{sub:first-order-unaligned}, and define $\N_{h,1}^z(\yy_0):=\{{{\mathbf{y}}}_{\Delta}(z)\}$. 
We also denote by  $\{\yy_{\Delta,i}(z)\}_{i=1}^4$ and $(\alpha_2(z),\beta_2(z))$ the four grid nodes surrounding $\yy_0(z)$ and the relative grid shift parameters respectively, as defined in Section~\ref{sub:second-order-correction}, and define $\N_{h,4}^z( \yy_0):=\{{{\mathbf{y}}}_{\Delta,i}(z)\}_{i=1}^4$. 

From the definition in \eqref{eq:svdefine} we can compute the expansion \eqref{eq:s-expansion} with $q=1$ and find
\begin{equation}\label{eq:expansion-s-3D}
    s({\mathbf{y}};z) = s_{0}({\mathbf{y}};z)+ s_{1}({\mathbf{y}};z)+\Oo(\vert{\mathbf{y}}\vert),
\end{equation}
with $s_{k}({\mathbf{y}};z)=\vert{\mathbf{y}}\vert^{k-1}\phi_{k}({\mathbf{y}}/\vert{\mathbf{y}}\vert;z)$,
$k=0,1$.
The expressions for $s_k$ are given in 
Theorem~\ref{thm:kernelexpansions} below.
 We can then apply the additive splitting \eqref{eq:additive-splitting}:
\begin{align*}
    \mathcal{U}^3_h[f(\cdot,\cdot,z)] =&\, h^2\sum_{{\mathbf{y}}\in h\mathbb{Z}^2 \setminus \N_{h,4}^z(\yy_0)}f({\mathbf{y}}, z) \nonumber\\ 
    & + h\sum_{i=1}^4 \omega_i[{s_{0}(\,\cdot\,;z)};\alpha_2(z),\beta_2(z)] v(\yy_{\Delta,i}(z),z) \nonumber\\ 
    & + h^2\, \omega[s_{1}(\,\cdot\,; z);\alpha_1(z),\beta_1(z)]\,v({{{\mathbf{y}}}}_{\Delta}(z), z)\,,\nonumber\\
    & + h\sum_{{\mathbf{y}}\in \N_{h,4}^z(\yy_0) \setminus \N_{h,1}^z(\yy_0) } \left\{ s({\mathbf{y}}-\yy_0(z);z)-{s_{0} ({\mathbf{y}}-\yy_0(z);z)}\right\}  v({\mathbf{y}}, z)   \nonumber.
\end{align*}
Then the three-dimensional third order method $\mathcal{V}^{3,z}_h$, obtained by applying $\mathcal{U}^3_h$ plane-by-plane along the $z$-direction, is given by
\begin{align}
     \mathcal{V}^{3,z}_h&[f] :=\, h\sum_{k\in\mathbb{Z}} \mathcal{U}^3_h[f(\cdot,\cdot,z_k)] = h^3\sum_{\bar{\mathbf{y}}\in h\mathbb{Z}^3\, \setminus\, \left(\bigcup_{k\in\mathbb{Z}} \N_{h,4}^{z_k}(\yy_{0}(z_k))\right)}f(\bar{\mathbf{y}})\nonumber\\
    &+ h^2\sum_{k\in\mathbb{Z}}\sum_{i=1}^4 \omega_i[{s_{0}(\,\cdot\,;z_k)};\alpha_2(z_k),\beta_2(z_k)] v(\yy_{\Delta,i}(z_k), z_k) \label{eq:Q3z-3D-general}\\
    &+ h^3\sum_{k\in\mathbb{Z}} \omega[s_{1}(\,\cdot\,;z_k);\alpha_1(z_k),\beta_1(z_k)]\,v(\yy_{\Delta}(z_k), z_k)\nonumber\\
    &+ h^2\sum_{k\in\mathbb{Z}}\,\sum_{{\mathbf{y}}\in\N_{h,4}^{z_k}(\yy_{0})\setminus \N_{h,1}^{z_k}(\yy_{0})} \Big\{ s({\mathbf{y}}-\yy_{0}(z_k);z_k) - s_{0}({\mathbf{y}}-\yy_{0}(z_k);z_k) \Big\}  v({\mathbf{y}}, z_k)  \nonumber.
\end{align}

If we apply the two-dimensional rule $\mathcal{U}^3_h$ plane-by-plane along the $x$- or $y$- direction we obtain the corresponding rules $\mathcal{V}_h^{3,x}$ and $\mathcal{V}_h^{3,y}$ respectively. These cases are discussed in the following Section~\ref{sub:change-coordinates}. 
 
\subsubsection{Plane-by-plane correction for different normal directions}\label{sub:change-coordinates}

The normal direction $\nn$ directly affects the
decomposition of a three dimensional Cartesian grid into union of planes, on which we apply the new correction. 
We identify the dominant direction of $\nn=(n_x,n_y,n_z)$, and discretize the volumetric integral along that direction. If the dominant direction is $n_z$, 
the setup is the one described above. If it is $n_y$, we discretize along the $y$-direction,
and if it is $n_x$, we discretize along the $x$-direction.

We shall use $\mathcal{V}_h^p$ for the general three-dimensional $p$-order corrected trapezoidal rule, and using the division presented for the change of coordinates, we define it as:
\begin{align}
    \mathcal{V}_h^p[f]:=&
    \begin{cases}
    \mathcal{V}_h^{p,z}[f]=h\sum_{k\in\mathbb{Z}}\mathcal{U}_h^p[f(\cdot,\cdot,z_k)], & \text{ if $n_z$ is dominant},\\[0.2cm]
    \mathcal{V}_h^{p,y}[f]=h\sum_{k\in\mathbb{Z}}\mathcal{U}_h^p[f(\cdot,y_k,\cdot)], & \text{ if $n_y$ is dominant},\\[0.2cm]
    \mathcal{V}_h^{p,x}[f]=h\sum_{k\in\mathbb{Z}}\mathcal{U}_h^p[f(x_k,\cdot,\cdot)], & \text{ if $n_x$ is dominant}.
    \end{cases}\label{eq:Q-3D-general}
\end{align}

\subsection{Expansions of layer kernels} 
\label{sub:layer-kernels-approx}

In this section we will analyze
and expand the
singular functions defined in
\eqref{eq:svdefine}
when $K$ are the Laplace kernels
\eqref{eq:layer-kernels}:
\begin{equation}\label{eq:layer-kernels-convenient-setting}
    \begin{array}{rll}
        \text{(SL)\,:}& 
        s^{SL}({\mathbf{y}};z) &= \dfrac{1}{4\pi}\dfrac{1}{\vert P_\Gamma 
        ({\mathbf{y}}+\yy_0(z),z)-\bar\xx^*\vert},\\[0.4cm]
        \text{(DL)\,:}& s^{DL}({\mathbf{y}};z) &= -\dfrac{1}{4\pi}\dfrac{\nn_y^T(
        P_\Gamma 
        ({\mathbf{y}}+\yy_0(z),z)-\bar\xx^*)}
        {\vert P_\Gamma 
        ({\mathbf{y}}+\yy_0(z),z)-\bar\xx^*\vert^3},\\[0.4cm]
        \text{(DLC)\,:}& s^{DLC}({\mathbf{y}};z) &= \dfrac{1}{4\pi}\dfrac{\nn_x^T(
        P_\Gamma 
        ({\mathbf{y}}+\yy_0(z),z)-\bar\xx^*)}
        {\vert P_\Gamma 
        ({\mathbf{y}}+\yy_0(z),z)-\bar\xx^*\vert^3}.
    \end{array}
\end{equation}

The approach developed in Section~\ref{sec:corrected-trapezoidal-rules} requires analytic formulae of the expansions. This means that in order to adopt
the third order quadrature rule \eqref{eq:Q-3D-general} for the implicit boundary integral defined in
\eqref{eq:volume_potential}, one needs explicit analytical 
expressions for the first two expansion functions in \eqref{eq:expansion-s-3D}
related to
the singular functions above.
Through a third order approximation of the surface near the target point $\bar\xx^*$ we find these functions, which are given
in the following theorem.

\begin{theorem}\label{thm:kernelexpansions}
Let $\bar\xx^*\in\Gamma$ be the target point.
Suppose that the normal $\nn$ at $\bar\xx^*$ satisfies $\nn^T \bar\ee_z\neq 0$,
and that $\bar\yy_0(z)\in T_\varepsilon$.
Then, there is an $r>0$, depending on $z$, such that
all the singular functions defined in \eqref{eq:layer-kernels-convenient-setting} can
be written
in the form
\begin{equation}\label{eq:ellform}
s^{X}( {\mathbf{y}};z)=\dfrac{1}{\vert{\mathbf{y}}\vert}\ell^X\left(\vert {\mathbf{y}}\vert,\dfrac{ {\mathbf{y}}}{\vert {\mathbf{y}}\vert};z\right), \quad \vert{\mathbf{y}}\vert<r,
\qquad \text{X=SL, DLC, DL},
\end{equation}
where $\ell^X\in C^\infty((-r,r)\times {\mathbb S}^1)$.
Moreover, the functions
$s_0^{X}( {\mathbf{y}};z)$ and $s_1^{X}( {\mathbf{y}};z)$
in the expansion \eqref{eq:expansion-s-3D} are
\begin{equation}\label{eq:s0s1-kernels}
\begin{array}{ll}
    s_0^{SL}({\mathbf{y}}) := \dfrac{1}{\vert{\mathbf{y}}\vert}\dfrac{1}{\psi_0(\hat{\mathbf{y}})} \,, &
    s_1^{SL}({\mathbf{y}}) := -\dfrac{\psi_1(\hat{\mathbf{y}})}{\psi_0(\hat{\mathbf{y}})^2} \,,\\[0.5cm]
    s_0^{DLC}({\mathbf{y}}) := \dfrac{1}{\vert{\mathbf{y}}\vert}\dfrac{\xi_0(\hat{\mathbf{y}})}{\psi_0(\hat{\mathbf{y}})^3} \,, &
    s_1^{DLC}({\mathbf{y}}) := -3\dfrac{\xi_0(\hat{\mathbf{y}})\psi_1(\hat{\mathbf{y}})}{\psi_0(\hat{\mathbf{y}})^4}+\dfrac{\xi_1(\hat{\mathbf{y}})}{\psi_0(\hat{\mathbf{y}})^3} \,,\\[0.5cm] %s_0^{DL}(\hat{\mathbf{y}}) &:= \dfrac{\xi_0(\hat{\mathbf{y}})}{\psi_0(\hat{\mathbf{y}})^3} \\
    s_0^{DL}({\mathbf{y}}) := s_0^{DLC}({\mathbf{y}}) \,, & 
    s_1^{DL}({\mathbf{y}}) := -3\dfrac{\xi_0(\hat{\mathbf{y}})\psi_1(\hat{\mathbf{y}})}{\psi_0(\hat{\mathbf{y}})^4}+\dfrac{\tilde \xi_1(\hat{\mathbf{y}})}{\psi_0(\hat{\mathbf{y}})^3}\,.
\end{array}
\end{equation}
where $\hat{\mathbf{y}}={\mathbf{y}}/\vert{\mathbf{y}}\vert$ and
$\psi_j$, $\xi_j$
and $\tilde\xi_1$ are given explicitly in
Section~\ref{subsub:kernel-expansions-summary}.
\end{theorem}

In order to use the expansions in the theorem in our
quadrature method, we need to be able to evaluate
the functions  $\psi_j$, $\xi_j$
and $\tilde\xi_1$. They depend on the 
local behavior of $\Gamma$ at the target point $\bar\xx^*$,
more precisely on the principal directions and curvatures, and the third derivatives of the function whose graph locally describes $\Gamma$.
In Appendix~\ref{sec:appendixB} 
it is described how those quantities can be
computed numerically using the closest point mapping.

In the subsequent subsections we will prove 
Theorem~\ref{thm:kernelexpansions}.
First, in Section~\ref{subsub:proj-mapping-to-layer-kernels}, we rotate the frame of reference and look at $\Gamma$ locally as the graph of a two-dimensional function.  
Second, we expand the expressions we obtained around ${\mathbf{y}}=0$ in
Section~\ref{subsub:fbh-expand}
and apply a general lemma to show \eqref{eq:ellform}
in Section~\ref{subsub:kernel-expandability}.
Finally, we use the expansions to 
derive expressions for
$s_0^{X}( {\mathbf{y}};z)$ and $s_1^{X}( {\mathbf{y}};z)$ in Section~\ref{subsub:kernel-expansions-summary}.

\subsubsection{Expressions of the layer kernels via the projection mapping}\label{subsub:proj-mapping-to-layer-kernels}

Let $\bar\xx^*\in\Gamma$ be the target point. At $\bar\xx^*$ we
denote the surface
principal directions  $\bar{\bm{\tau}}_1,\bar{\bm{\tau}}_2$,
the normal $\nn$, and 
the principal curvatures $\ka_1,\ka_2$. 
We introduce the
\emph{principal basis} 
$B=(\bar{\bm{\tau}}_1,\bar{\bm{\tau}}_2,\nn)$ and the notation
\[
    (x'_1, x'_2, x_3')_B
    := \bar{\bm{\tau}}_1 x_1' + \bar{\bm{\tau}}_2 x_2' + \nn x_3'.
\]
The basis vectors used here are assumed to be normalized.
If $\bar\xx'$
are the coordinates in the $B$-basis
for the point $\bar\xx$ in the canonical basis
$(\bar\ee_x,\bar\ee_y,\bar\ee_z)$, we denote by
$Q$ the (orthogonal) change of basis matrix, satisfying
\begin{equation}\label{eq:def-Q}
 \bar\xx = Q\bar\xx',
 \qquad Q = 
     \left(
    \begin{array}{ccc}
        \vert & \vert & \vert \\ \bar{\bm{\tau}}_1 & 
        \bar{\bm{\tau}}_2& \nn \\
        \vert & \vert & \vert
    \end{array}
    \right), \qquad Q^TQ=I.
\end{equation}
The surface $\Gamma$
can now be parameterized locally in
the $B$-coordinates. More precisely,
 in a neighborhood of the origin,
${\mathcal I}_{L'}=\{{\mathbf{y}}'\in{\mathbb{R}}^2\ \vert\ 
\vert{\mathbf{y}}'\vert< L'\}$, we can represent
$\Gamma$ as the image of a smooth
function $f:{\mathbb{R}}^2\to{\mathbb{R}}$ with $f\in C^\infty({\mathcal I}_{L'})$ such that
$$
\bar\xx^*+({\mathbf{y}}',f({\mathbf{y}}'))_B \in \Gamma.
$$
The constant $L'$ 
depends on the maximum
curvature of $\Gamma$ and
can be taken to be independent
of $\bar\xx^*$ .
Moreover,
since $P_\Gamma$ is smooth in the
tubular neighborhood $T_\varepsilon$
of $\Gamma$, the mapping
$({\mathbf{y}}',z')\mapsto P_\Gamma(\bar\xx^*+({\mathbf{y}}',z')_B)$
is smooth for
$({\mathbf{y}}',z')\in {\mathcal T}_{\varepsilon}
:= \{({\mathbf{y}}',z') \in{\mathbb{R}}^3\,:\, \bar\xx^*+({\mathbf{y}}',z')_B\in T_\varepsilon\}$.
Therefore, for
$({\mathbf{y}}',z')\in {\mathcal M}_{L}={\mathcal T}_{\varepsilon}
\cap ({\mathcal I}_{L}\times{\mathbb{R}})$,
with $L$ possibly smaller than $L'$,
we can use the $B$-basis
and $f$,
to write
the closest point mapping
as
\begin{equation}\label{eq:Pgamma-to-f-and-h}
P_\Gamma(\bar\xx^*+({\mathbf{y}}',z')_B) = \bar\xx^*+ \big( \yy_{\text{p}},f(\yy_{\text{p}})\big)_B,
\qquad \yy_{\text{p}} := \bh({\mathbf{y}}',z'),
\end{equation}
for some smooth function $\bh\in C^\infty({\mathcal M}_{L})$.
The constant $L$ is chosen such that
\[
\sup_{({\mathbf{y}}',z')\in {\mathcal M}_{L}} \vert\bh({\mathbf{y}}',z')\vert\leq L'.
\]
Clearly $\bh(\zero,z')=\zero$,
which guarantees that $L>0$.

We now
write $({\mathbf{y}}+\yy_0(z),z)$ as a point in the $B$-basis centered in
the target point $\bar\xx^*$,
\[
({\mathbf{y}}+\yy_0(z),z) = \bar\xx^* + ({\mathbf{y}}',z')_B.
\]
For $({\mathbf{y}}',z')\in{\mathcal M}_{L}$
we can then write the numerators and denominators of the layer kernels \eqref{eq:layer-kernels-convenient-setting} using \eqref{eq:Pgamma-to-f-and-h}
and the orthogonality
of $Q$:
\begin{align}
    \left\vert P_\Gamma({\mathbf{y}}+\yy_0(z),z)-\bar\xx^*\right\vert =& \left\vert \big( \yy_{\text{p}},f(\yy_{\text{p}}) \big)_B \right\vert = \left\vert \big( \yy_{\text{p}},f(\yy_{\text{p}}) \big) \right\vert, \label{eq:SL-denominator-explicit}\\[0.3cm]
    \nn_x^T( P_\Gamma ({\mathbf{y}}+\yy_0(z),z)-\bar\xx^*) =& 
    \left(
    \begin{array}{c}
        \zero \\ 1
    \end{array}
    \right)^T_B
    \left(
    \begin{array}{c}
        \yy_{\text{p}} \\ f(\yy_{\text{p}})
    \end{array}
    \right)_B = \left(
    \begin{array}{c}
        \zero \\ 1
    \end{array}
    \right)^T \left(
    \begin{array}{c}
        \yy_{\text{p}} \\ f(\yy_{\text{p}})
    \end{array}
    \right),\label{eq:DLC-numerator-explicit}\\[0.4cm]
\nn_y^T( P_\Gamma ({\mathbf{y}}+\yy_0(z),z)-\bar\xx^*) =& \frac{1}{\sqrt{1+(\nabla f(\yy_{\text{p}}))^2}}
    \left(
    \begin{array}{c}
        -\nabla f(\yy_{\text{p}}) \\ 1
    \end{array}
    \right)^T_B
    \left(
    \begin{array}{c}
        \yy_{\text{p}} \\ f(\yy_{\text{p}})
    \end{array}
    \right)_B \nonumber\\
    =& \frac{1}{\sqrt{1+(\nabla f(\yy_{\text{p}}))^2}}
    \left(
    \begin{array}{c}
        -\nabla f(\yy_{\text{p}}) \\ 1
    \end{array}
    \right)^T
    \left(
    \begin{array}{c}
        \yy_{\text{p}} \\ f(\yy_{\text{p}})
    \end{array}
    \right).\label{eq:DL-numerator-explicit}
\end{align}

We next have to find how $({\mathbf{y}}',z')$ depends on 
${\mathbf{y}}$ and $z$.
From the definitions above we have
$$
   \left(
    \begin{array}{c}
        {\mathbf{y}}\\ 0
    \end{array}
    \right)
    +\bar\yy_0(z)=
   \left(
    \begin{array}{c}
        {\mathbf{y}}+\yy_0(z) \\ z
    \end{array}
    \right)
 = \bar\xx^* + Q
   \left(
    \begin{array}{c}
        {\mathbf{y}}' \\ z'
    \end{array}
    \right).
$$
Since $\bar\yy_0(z)-\bar\xx^*$
is parallel to the normal $\nn$
by definition,
we can express this as
$$
\bar\yy_0(z)-\bar\xx^* = (0,\eta(z))_B, \quad \Rightarrow \quad 
\left(
\begin{array}{c}
    {\mathbf{y}}\\ 0
\end{array}
\right)
= Q
\left(
\begin{array}{c}
    {\mathbf{y}}' \\ z'-\eta(z)
\end{array}
\right).
$$
where 
$\eta(z):=d_\Gamma(\bar\yy_0(z))$ is the signed distance of $\bar\yy_0(z)$ to $\Gamma$.
Defining
\begin{equation}\label{eq:def-QT-A-dvec}
Q^{T}=\left(
    \begin{array}{cc}
        {A} & \mathbf{c} \\ \mathbf{d}^T & \alpha
    \end{array}
    \right),\qquad
    {A} \in{\mathbb{R}}^{2\times 2},\quad \mathbf{c},\mathbf{d}\in{\mathbb{R}}^{2\times 1}, \quad \alpha\in{\mathbb{R}},
\end{equation}
we finally obtain
\begin{equation}\label{eq:x=Ay and z}
\left\{
\begin{array}{rl}
    {\mathbf{y}}' =& {A}{\mathbf{y}}\,, \\
    z'  =& \mathbf{d}^T{\mathbf{y}}+\eta(z).
\end{array}
\right.
\end{equation}
Therefore, we can write the kernels \eqref{eq:layer-kernels-convenient-setting} using (\ref{eq:SL-denominator-explicit},\ref{eq:DLC-numerator-explicit},\ref{eq:DL-numerator-explicit}) and \eqref{eq:x=Ay and z}:
\begin{equation}\label{eq:layer-kernels-not-expanded}
    \begin{array}{lrl}
        & \yy_{\text{p}}:=& \bh(A{\mathbf{y}},\,\mathbf{d}^T{\mathbf{y}}+\eta(z)), \\[0.2cm]
        \text{(SL) }& s^{SL}({\mathbf{y}};z) = & \dfrac{1}{4\pi}\dfrac{1}{\vert(\yy_{\text{p}},f(\yy_{\text{p}}))\vert},\\[0.5cm]
        \text{(DLC) }& s^{DLC}({\mathbf{y}};z) = & \dfrac{1}{4\pi}\dfrac{f(\yy_{\text{p}})}{\vert(\yy_{\text{p}},f(\yy_{\text{p}}))\vert^3}, \\[0.5cm]
        \text{(DL) }& s^{DL}({\mathbf{y}};z) = & -\dfrac{1}{4\pi}\dfrac{(-\nabla f(\yy_{\text{p}}),1)}{\vert(\yy_{\text{p}},f(\yy_{\text{p}}))\vert^3\sqrt{1+(\nabla f(\yy_{\text{p}}))^2}}
        \left( \begin{array}{c} \yy_{\text{p}} \\ f(\yy_{\text{p}}) \end{array} \right) . 
    \end{array}
\end{equation}
These expressions are
valid for $({\mathbf{y}}',z')\in{\mathcal M}_L$.
By \eqref{eq:x=Ay and z} 
and the fact that $\vert A{\mathbf{y}}\vert\leq \vert Q^T({\mathbf{y}},0)\vert=\vert{\mathbf{y}}\vert$
(with equality if ${\mathbf{y}}\perp\mathbf{d}$)
they are 
therefore valid
when $({\mathbf{y}}+\yy_0(z),z)\in T_\varepsilon$ and 
$\vert{\mathbf{y}}\vert< L$.

\subsubsection{Expansion of $f$ and $\bh$}\label{subsub:fbh-expand}

By the definition of the $B$-basis,
the function $f$ introduced above in
Section~\ref{subsub:proj-mapping-to-layer-kernels}
satisfies
\begin{equation}\label{eq:f-properties-def-M}
f(\zero)=0,\ \qquad \nabla f(\zero)=\zero,\ \qquad \dfrac{\partial^2 f}{\partial \xx^2}(\zero)=
\left(
    \begin{matrix}
    \ka_1 & 0 \\ 0 & \ka_2
    \end{matrix}
\right)=:M.
\end{equation}
The Taylor expansions up to second order for $f$
and $\nabla f$ are then given by
\begin{equation}\label{eq:high_order_f_df}
    \begin{array}{rl}
    f({\mathbf{y}}) &= \frac{1}{2}{\mathbf{y}}^TM{\mathbf{y}}+B({\mathbf{y}},{\mathbf{y}},{\mathbf{y}})+\Oo(\vert{\mathbf{y}}\vert^4)\,, \\[0.1cm]
    \nabla f({\mathbf{y}}) &= M{\mathbf{y}}+C({\mathbf{y}},{\mathbf{y}})+\Oo(\vert{\mathbf{y}}\vert^3),
    \end{array}
\end{equation}
where $B$ is the third order trilinear term, and $C$ is its bilinear gradient.
With ${\mathbf{y}}=(x,y)$, they are given by
\begin{equation}\label{eq:surface-B-C}
\begin{array}{rl}
    B({\mathbf{y}},{\mathbf{y}},{\mathbf{y}})&%\equiv B(x,y) 
    :=\frac{1}{2}\left[ f_{xxx}\frac{x^3}{3}+f_{yyy}\frac{y^3}{3}+f_{xxy}x^2y+f_{xyy}xy^2 \right],\\[0.3cm]
    C({\mathbf{y}},{\mathbf{y}})& %\equiv C(x,y) 
    :=
    \left(
    \begin{matrix}
    \frac{\partial B}{\partial x} \\[0.2cm] 
    \frac{\partial B}{\partial y}
    \end{matrix}
    \right)=\dfrac{1}{2}
    \left(
    \begin{matrix}
    f_{xxx}x^2+2f_{xxy}xy+f_{xyy}y^2 \\[0.2cm] f_{yyy}y^2+2f_{xyy}xy+f_{xxy}x^2
    \end{matrix}
    \right);
\end{array}
\end{equation}
where $f_{xxx},f_{xxy},f_{xyy},f_{yyy}$ are the third order derivatives of $f$ evaluated in $\zero$.

We next need to expand $\bh$.
It is given by the following lemma,
the proof of which can be found in the
Appendix~\ref{sub:A:properties-of-h}.
\begin{lemma}\label{lem:properties-of-h}
Let
\[
D(z')=(I-z'M)^{-1},
\]
with $M$ given in \eqref{eq:f-properties-def-M}.
For $({\mathbf{y}}',z')\in {\mathcal M}_L$ the matrix
is well-defined.
The
function $\bh\in C^{\infty}(
{\mathcal M}_L)$
introduced in Section~\ref{subsub:proj-mapping-to-layer-kernels}
then satisfies
\begin{equation}\label{eq:hprop1}
\bh(\zero,z')=\zero, \qquad 
\frac{\partial \bh}{\partial z}(\zero,z')=\zero, \qquad
\frac{\partial \bh}{\partial {\mathbf{y}}}(\zero,z')=D(z'), 
\end{equation}
and the Taylor expansion of $\bh$ can be written in the form
\begin{equation}%\label{eq:hprop2}
\bh({\mathbf{y}}',z') \,=\, D(z'){\mathbf{y}}' + z'\,D(z')C\Bigl(D(z'){\mathbf{y}}',D(z'){\mathbf{y}}'\Bigr) +\Oo(\vert{\mathbf{y}}'\vert^3),\label{eq:hig_order_proj_x}
\end{equation}
where $C$ is defined in \eqref{eq:surface-B-C}.
\end{lemma}

\subsubsection{General form of the kernels}\label{subsub:kernel-expandability}

We now have expressions \eqref{eq:layer-kernels-not-expanded} of the kernels and expansions around ${\mathbf{y}}=\yp=\zero$
of $\bh$ and $f$.
The next step is to prove \eqref{eq:ellform}, i.e. that
the three kernels in
\eqref{eq:layer-kernels-not-expanded}
can all be written in the form $\vert{\mathbf{y}}\vert^{-1}\ell(\vert{\mathbf{y}}\vert,{\mathbf{y}}/\vert{\mathbf{y}}\vert)$. 
To do this we use the following lemma, a proof of which can be found in Appendix~\ref{sub:A:g=f-over-distance}.

\begin{lemma}\label{lem:g=f-over-distance}
Let $\bar\bg:{\mathbb{R}}^m\to{\mathbb{R}}^n$ be $C^\infty(B_{r_0}(\zero))$ for some $r_0>0$, with
$n>m$, $\bar\bg(\zero)=\zero$,
and $D \bar\bg(\zero)\in{\mathbb{R}}^{n\times m}$ has full rank. Let $\bar\bp:{\mathbb{R}}^m\to{\mathbb{R}}^n$ be $C^\infty(B_{r_0}(\zero))$, such that $\bar\bp(\zero)^TD\bar\bg(\zero)=\zero$. Then there exist functions $\ell_1,\,\ell_2$ 
and $0<r_1\leq r_0$
such that $\ell_i:{\mathbb{R}}\times\s^{m-1}\to{\mathbb{R}}$, $\ell_i\in C^\infty((-r_1,r_1)\times\s^{m-1})$, $i=1,2$ and
\[
\dfrac{1}{\vert\bar\bg({\mathbf{y}})\vert}=\dfrac{1}{\vert{\mathbf{y}}\vert}\ell_1\left(\vert{\mathbf{y}}\vert,\dfrac{{\mathbf{y}}}{\vert{\mathbf{y}}\vert}\right) , \qquad \dfrac{\bar\bp({\mathbf{y}})^T\bar\bg({\mathbf{y}})}{\vert\bar\bg({\mathbf{y}})\vert^3}=\dfrac{1}{\vert{\mathbf{y}}\vert}\ell_2\left(\vert{\mathbf{y}}\vert,\dfrac{{\mathbf{y}}}{\vert{\mathbf{y}}\vert}\right).
\]
\end{lemma}

For the single-layer kernel we take 
\[
\bar\bg({\mathbf{y}})=(\yp,f(\yp))
=\Big( \bh(A{\mathbf{y}},\mathbf{d}^T{\mathbf{y}}+\eta(z)),f\big(\bh(A{\mathbf{y}},\mathbf{d}^T{\mathbf{y}}+\eta(z))\big) \Big).
\] 
For $(\zero,\eta(z))\in\mathcal{M}_L$,
i.e. when $\vert\eta(z)\vert<\varepsilon$,
Lemma~\ref{lem:properties-of-h} 
gives that
$\bar\bg(\zero)=(\zero,0)$ and
\begin{align*}
&\dfrac{\partial \bar\bg}{\partial {\mathbf{y}}}(\zero) = 
\left(
\begin{array}{cc}
     \frac{\partial \bh}{\partial {\mathbf{y}}}(\zero,\eta(z))
     \left(A + \frac{\partial \bh}{\partial z}(\zero,\eta(z))\mathbf{d}^T\right)\\
     \left(\frac{\partial \bh}{\partial {\mathbf{y}}}(\zero,\eta(z))\left(A + \frac{\partial \bh}{\partial z}(\zero,\eta(z))\mathbf{d}^T\right)\right)^T\nabla f(\zero)  
\end{array}
\right)
= 
\left(
\begin{array}{cc}
     D(\eta(z))A  \\
     \zero  
\end{array}
\right),
\end{align*}
which has full rank
since $\det A = \bar\ee_z^T(
\bar{\bm{\tau}}_1\times\bar{\bm{\tau}}_2)
=\bar\ee_z^T\nn\neq 0$.
Hence, 
\eqref{eq:layer-kernels-not-expanded} together with
the first result of Lemma~\ref{lem:g=f-over-distance}
now shows \eqref{eq:ellform} for $X=SL$.

For the double-layer case, we 
let $\bar\bp({\mathbf{y}})=(-\nabla f(\yp),1)/\sqrt{1+\vert\nabla f(\yp)\vert^2}$
so that $s^{DL}=-\bar\bp^T\bar\bg/4\pi\vert\bar\bg\vert^3$ by
\eqref{eq:layer-kernels-not-expanded}.
Then Lemma~\ref{lem:properties-of-h} gives
\[
\bar\bp(\zero)^TD\bar\bg(\zero) = \bar\bp(\zero)^T\dfrac{\partial \bar\bg}{\partial {\mathbf{y}}}(\zero)  = 
\left(
\begin{array}{c}
     -\nabla f(\zero)  \\
     1  
\end{array}
\right)^T
\left(
\begin{array}{cc}
     D(\eta(z))A  \\
     \zero  
\end{array}
\right)
= \left(
\begin{array}{c}
     \zero  \\
     0  
\end{array}
\right),
\]
and the second result of Lemma \ref{lem:g=f-over-distance}
shows \eqref{eq:ellform} for $X=DL$.
Finally, for the double layer conjugate kernel we take simply
$\bar\bp({\mathbf{y}})=(0,0,1)$, which again makes $\bar\bp(\zero)^TD\bar\bg(\zero) =\zero$
and \eqref{eq:ellform} for $X=DLC$ follows as before.
This completes the proof of \eqref{eq:ellform}.

\subsubsection{Kernel expansions}\label{subsub:kernel-expansions-summary}
The expansion of the kernels is based on the expansions of
$f$ in \eqref{eq:high_order_f_df} and $\bh$ in \eqref{eq:hig_order_proj_x}.
We will skip most tedious intermediate calculations and focus on the
end results.

We recall that
$$
        {\mathbf{y}}'={A}{\mathbf{y}},\qquad z'=\mathbf{d}^T{\mathbf{y}}+\eta(z),\qquad \eta(z)=d_\Gamma(\bar\yy_0(z)),
        \qquad \hat{\mathbf{y}} := {\mathbf{y}}/\vert{\mathbf{y}}\vert.
$$
In the first step we expand the functions $f(\bh({\mathbf{y}}',z'))$, $D(z')$ and $\bh({\mathbf{y}}',z')$ as functions
of ${\mathbf{y}}$, instead of $\yp$ and ${\mathbf{y}}'$ as before. We get
\begin{align*}
        D(z') &= D_0\left[ I+\mathbf{d}^T{\mathbf{y}} D_0M\right]+\Oo(\vert{\mathbf{y}}\vert^2),\\
\bh({\mathbf{y}}',z') 
        &= \chi_0(\hat {\mathbf{y}})\vert{\mathbf{y}}\vert+\chi_1(\hat {\mathbf{y}})\vert{\mathbf{y}}\vert^2+\Oo(\vert{\mathbf{y}}\vert^3),\\
        f(\yy_{\text{p}}) 
        &= \xi_0(\hat{\mathbf{y}})\vert{\mathbf{y}}\vert^2+\xi_1(\hat{\mathbf{y}})\vert{\mathbf{y}}\vert^3+\Oo(\vert{\mathbf{y}}\vert^4),
\end{align*}
where $D_0:=(I-\eta M)^{-1}$,
\begin{equation*}
    \begin{array}{rl}
            \chi_0 ( {\mathbf{y}}):=& D_0A {\mathbf{y}}, \\[0.2cm]
            \chi_1 ( {\mathbf{y}}):=& (\mathbf{d}^T{\mathbf{y}})D_0D_0MA{\mathbf{y}} +\eta D_0 C(D_0A{\mathbf{y}},D_0A{\mathbf{y}}),\\[0.2cm]
            \xi_0({\mathbf{y}}):=& \dfrac{1}{2}{\mathbf{y}}^T(A^TD_0^TMD_0A){\mathbf{y}}, \\[0.2cm]
            \xi_1({\mathbf{y}}):=& \dfrac{1}{2}\eta(D_0C(D_0A{\mathbf{y}},D_0A{\mathbf{y}}))^TMD_0A{\mathbf{y}}
            + (\mathbf{d}^T{\mathbf{y}}){\mathbf{y}}^TA^T(M^TD_0^TD_0^TMD_0)A{\mathbf{y}}\\ 
            & + \dfrac{1}{2}\eta(D_0A{\mathbf{y}})^TMD_0C(D_0A{\mathbf{y}},D_0A{\mathbf{y}}) + B(D_0A{\mathbf{y}},D_0A{\mathbf{y}},D_0A{\mathbf{y}}).
    \end{array}
\end{equation*}
In this step we used the fact that $\chi_j$ and $\xi_j$ are
homogeneous of degree $j+1$ and $j+2$ respectively, 
so that $\chi_j({\mathbf{y}})=\chi_j(\hat{\mathbf{y}})\vert{\mathbf{y}}\vert^{j+1}$
and $\xi_j({\mathbf{y}})=\xi_j(\hat{\mathbf{y}})\vert{\mathbf{y}}\vert^{j+2}$.
From these expansions for $f$ and $\bh$ we obtain furthermore that
    \begin{align}\label{eq:Pexpansions}
        \big\vert(\yy_{\text{p}},f(\yy_{\text{p}}))\big\vert &=
        \psi_0(\hat {\mathbf{y}})\vert{\mathbf{y}}\vert+\psi_1(\hat{\mathbf{y}})\vert{\mathbf{y}}\vert^2+\Oo(\vert{\mathbf{y}}\vert^3),\\
        \dfrac{(\nabla f(\yy_{\text{p}}),-1)}{\sqrt{1+(\nabla f(\yy_{\text{p}}))^2}} \left( \begin{array}{c}
        \yy_{\text{p}} \\ f(\yy_{\text{p}}) \end{array}\right) &= 
         \xi_0(\hat {\mathbf{y}})\vert{\mathbf{y}}\vert^2+\tilde \xi_1(\hat{\mathbf{y}})\vert{\mathbf{y}}\vert^3+\Oo(\vert{\mathbf{y}}\vert^4),
         \nonumber
    \end{align}
where
\begin{align*}
    \psi_0({\mathbf{y}}):=& \vert\chi_0( {\mathbf{y}})\vert,\qquad \psi_1({\mathbf{y}}):=\ \ \dfrac{\chi_0({\mathbf{y}})^T \chi_1({\mathbf{y}})}{\vert\chi_0({\mathbf{y}})\vert},\\
    \tilde \xi_1( {\mathbf{y}}) :=&  \frac{1}{2}\eta (D_0C(D_0A {\mathbf{y}},D_0A{\mathbf{y}}))^TMD_0A {\mathbf{y}} + (\mathbf{d}^T {\mathbf{y}}) {\mathbf{y}}^TA^TD_0MD_0D_0MA {\mathbf{y}} \\
    &- B(D_0A {\mathbf{y}},D_0A{\mathbf{y}},D_0A{\mathbf{y}}) + {\mathbf{y}}^TA^TD_0^T(I+\eta MD_0)C(D_0A {\mathbf{y}},D_0A{\mathbf{y}}) \\
    &- \frac{1}{2}\eta (D_0A {\mathbf{y}})^TMD_0C(D_0A {\mathbf{y}},D_0A{\mathbf{y}}) .
\end{align*}
Here $\tilde\xi_1$ is homogeneous of degree three. 
Using \eqref{eq:Pexpansions}
one can finally deduce the expansions of the kernels in \eqref{eq:layer-kernels-not-expanded}.
This concludes the proof of Theorem~\ref{thm:kernelexpansions}.

We note that the matrix $A$ and vector $\mathbf{d}$ contain elements
of the principal directions and normal at the target point;
see \eqref{eq:def-Q} and \eqref{eq:def-QT-A-dvec}.
The matrices $D_0$ and $M$
are built from the principal curvatures of $\Gamma$, 
and the functions $B$ and $C$ contain the
third derivatives of $f$;
see \eqref{eq:high_order_f_df} and \eqref{eq:surface-B-C}.
In Appendix~\ref{sec:appendixB} we show how to numerically compute the information about the surface in the target point 
($\kappa_1,\kappa_2$, $\bar{\bm{\tau}}_1$,$\bar{\bm{\tau}}_2$, and the third derivatives of $f$,  $f_{xxx}$, $f_{xxy}$, $f_{xyy}$, $f_{yyy}$) 
using the projection mapping $P_\Gamma$ and its derivatives.

\subsection{{Requirements for order higher quadratures for the singular IBIM integrals}}

{Given the class of singular integrands, the main obstruction to obtaining higher order quadratures using the proposed approach is the smoothness of the surface. When applying the proposed method in the IBIM formulation using uniform Cartesian grids, one needs firstly a sufficiently accurate approximation of the distance function to the surface, $d_\Gamma$, or the projection, $P_\Gamma$, on the grid nodes.}

{The construction of these functions are application dependent, but general methodologies do exist, see e.g. \cite{osher2006level}. 
If the surfaces are reconstructed on a grid by a level set method, then typically one does not expect that $d_\Gamma$ be more than 4th order accurate in the grid spacing due to the limitation imposed by commonly used level set reinitialization algorithms \cite{osher2006level}. This may cause a main bottleneck in practice.
Then one needs to extract the surface's geometrical information from finite differences of $d_\Gamma$ or $P_\Gamma$ -- in this paper, the related quantities to be approximated are the partial derivatives of $f$ defined in \eqref{eq:surface-B-C}, where $f$ is defined in \eqref{eq:kernelfunction}.
In Appendix~\ref{sec:appendixB} the reader will find more details.}

{When the surface is sufficiently smooth, it has a non-zero reach; i.e. $d_\Gamma$ is smooth within $T_{\tau_\Gamma}$ for some $\tau_\Gamma>0$. The Cartesian grid inside $T_{\tau_\Gamma}$ should be sufficiently dense to support the finite difference stencil around any node inside $T_\varepsilon$, where $\varepsilon\leq \tau_\Gamma$. 
Thus higher order approximations require denser grids around the surface to support the wider finite difference stencils used in high order finite differences.} {For example the second order corrected rule $\mathcal{V}_h^2$ needs curvature information which is obtained through a centered 5 point three-dimensional stencil. This implies that one needs accurate $d_\Gamma$ or $P_\Gamma$ within the distance of $\varepsilon+2h$ to the surface, and that $\varepsilon+2h$ should be smaller than the reach $\tau_\Gamma$.
Analogously, the third order information about the surface needed for $\mathcal{V}_h^3$ is obtained using a $5\times5\times5$ stencil around each node, which leads to the bound $\varepsilon+2\sqrt{2}h<\tau_\Gamma$. }
{If the surface geometry varies ``wildly", we envision that the proposed method should/could be generalized to 
multi-resolution gridding for efficiency. }  

{ We present an example supporting the above discussion in Section~\ref{sec:numerical-tests}. We also refer the interested readers to the results and discussion in the recent paper \cite{izzo2022PB}, for an application of the proposed quadratures in computing the electrostatic potentials of large molecules in a solvent. }

\section{Numerical tests}
\label{sec:numerical-tests}
In this Section we test the corrected trapezoidal rules derived in Section~\ref{sec:corrected-trapezoidal-rules} and Section~\ref{sec:application-layer-potentials-3D}. In Section~\ref{sec:num_2D} we test the rules $Q_h^p$ for integrating functions of the kind $s_k\,v$ from Section~\ref{sub:high-order-correction}, and then the general rules  $\,\mathcal{U}_h^p$ for integrating $s\,v$ from Section~\ref{sub:additive-splitting}. In Section~\ref{sub:num_3D} we test the third-order accurate quadrature rule $\mathcal{V}_h^3$ derived for the three-dimensional layer potentials discussed in Section~\ref{sec:application-layer-potentials-3D}.

\subsection{Corrections to the punctured trapezoidal rules in two dimensions}
\label{sec:num_2D}
The quadrature rules discussed in Section~\ref{sec:corrected-trapezoidal-rules} have been developed to correct any function of the kind 
\[
f(\xx) = s_k(\xx)v(\xx)\ ,\ \ s_k(\xx)=\vert\xx\vert^{k-1}\phi_k(\xx/\vert\xx\vert)\ ,\ \  k\in\mathbb{N}\setminus \{0\}\,,
\]
where $v$ is a smooth function, and then composite rules have been constructed to correct functions which can be expanded as \begin{align*}
f(\xx)&=s(\xx-\xx_0)v(\xx) \\
\text{where }\,s(\xx)&= s_0(\xx)+s_1(\xx)+s_2(\xx)+\dots.
\end{align*}
We tested the rules $Q_h^p$ for $p=1,2,3,4$ ($p=1$ \eqref{eq:single_correction_quadrature}, $p=2$ \eqref{eq:Q2-correction}, $p$ general \eqref{eq:Qp-correction}) for functions $s_k$, $k=0,1,2$. Specifically, we used the test function where $s_k$ and $v$ are:
\begin{align}\label{eq:test-sk}
    s_k(\xx) =& \vert\xx\vert^{k-1} \phi(\xx/\vert\xx\vert),\\[0.1cm]
    \phi(\xx/\vert\xx\vert) =& \phi(\cos(\psi(\xx)),\sin(\psi(\xx))) \nonumber\\[0.1cm]
    =& 4.2398+0.816735\cos(\psi(\xx)-0.2)-1.24397865\sin(2\psi(\xx)+0.1)\,, \nonumber\\
    v(\xx) =& \left(1.1+\Re\left(H_{\vert\xx\vert^2 +1}^{(1)}(3)\right)\right)\exp\left(-\vert\xx-(0.027,\,0.0197)\vert^8\right)\nonumber\\
    &  \cdot(0.5+\sin(\xx_1(\xx_2-1))).\nonumber
\end{align}
The function $H^{(1)}_\alpha$ is the Hankel function of the first kind of degree $\alpha$, and $\Re$ indicates the real part of a complex number. {
Although formally $v$ is not 
compactly supported, it is smaller than the
numerical machine precision outside $[-2,2]^2$, which we
use as integration domain.}

\begin{figure}
    \begin{center}
    \includegraphics[scale=1]{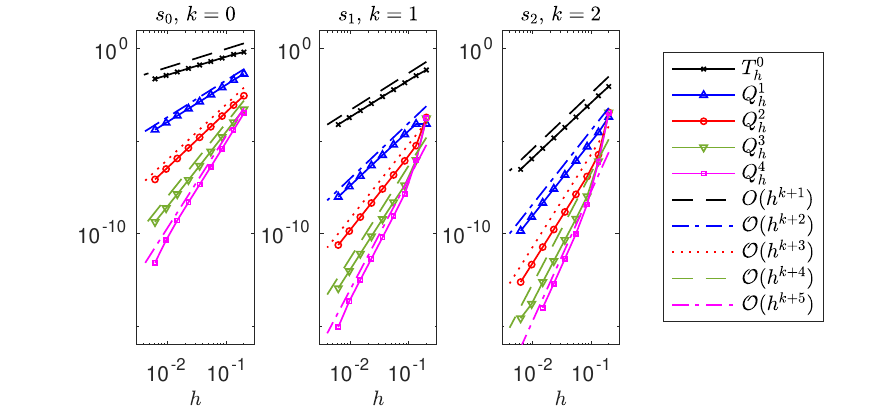}
    \end{center}
    \caption{\textbf{Correction of $s_k$ in two dimensions}. Error from integrating $s_k$ \eqref{eq:test-sk} with $p$-order correction $Q_h^p$. For $k=0,1,2$ (left, center, and right figures respectively) we present the difference
    between values obtained from 
    grid sizes $h$ and $h/1.5$,
    with the different methods. As expected, the order of accuracy is $k+p+1$ where $p$ is the order of the correction.}
    \label{fig:conv2D-sk}
\end{figure}

In Figure~\ref{fig:conv2D-sk} we plot the difference {between approximation values for grid sizes $h$
and $h/1.5$}, 
obtained for the 
four different quadratures $Q_h^p$, $p=1,2,3,4$, and the punctured trapezoidal rule $T_h^0$. The order of accuracy shown for integrating $s_k\,v$, $k=0,1,2$, is 
$k+1$ for the punctured trapezoidal rule and $k+p+1$ for the quadrature $Q_h^p$, as expected. The error constant is determined by the value of $(\alpha,\beta)$, and in our tests we fixed  $(\alpha,\beta)=(0.81,0.46)$.
The stencils used for the different quadratures are represented in Figure~\ref{fig:all_correction_grid_plot}. 
{The weights for the quadratures are all non-negative. Their
maximum values are shown in Table~\ref{table:weights}.
They are of moderate size also for the
high order corrections.}

\begin{table}
    \caption{\textbf{Maximum of the weights}. 
    {Largest weight $\max_{i=1}^{\tilde p} \omega_i$ in the stencil $\mathcal{N}_{h,\tilde p}$ for different correction orders $p=1,2,3,4$, and different singularity order $k=0,1,2$. The weights  correspond to the ones used in the tests shown in Figure~\ref{fig:conv2D-sk}. All weights are non-negative.}}
    \label{tab:weights-module}
    \begin{center}
    \begin{tabular}{c|c|c|c|c|}
        & \ $Q_h^1$\ & \ $Q_h^2$\ & \ $Q_h^3$\ & \ $Q_h^4$\ \\[0.1cm]
        \hline
        $k=0$ & 15.20855 & 11.39144 & 11.82856 & 11.61144  \\
        $k=1$ & 5.05848 & 4.91377 & 4.92476 & 5.11844 \\
        $k=2$ & 2.46476 & 4.59018 & 6.76066 & 8.88673 
    \end{tabular}
    \end{center}
    \label{table:weights}
\end{table}

\begin{figure}
    \begin{center}
        \includegraphics[scale=1]{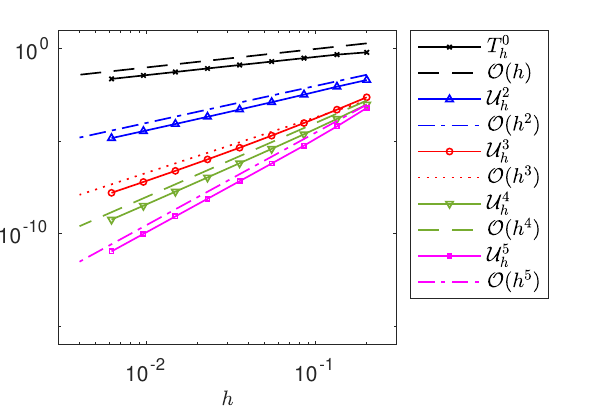}
    \end{center}
    \caption{\textbf{Corrected trapezoidal rules for a general function $s$ in two dimensions}. Corrected trapezoidal rules $\mathcal{U}_h^p$ for $p=2,3,4,5$ using additive splitting \eqref{eq:Qp-general-explicit} for the function $f=s\,v$ with singular integrand $s$ \eqref{eq:test-s-gen}. The first $p-1$ terms of the expansion \eqref{eq:s-expansion} ($s_k$, $k=0,1,\dots,p-2$) are needed to use $\mathcal{U}_h^p$. In the plot we see that the punctured trapezoidal rule $T_h^0$ has first order accuracy, and the corrections $\mathcal{U}_h^p$ have order of accuracy $p$ as predicted.}
    \label{fig:conv2D-s-gen}
\end{figure}

In order to test the general quadrature rule \eqref{eq:Qp-general-explicit} we used the function
\begin{align}
    s(\xx) =& \left(\vert\xx\vert^{-1}\phi_0(\xx)+\phi_1(\xx)+\vert\xx\vert\phi_2(\xx)+\vert\xx\vert^2\phi_3(\xx)+\vert\xx\vert^3r(\xx)\right)\label{eq:test-s-gen}\\[0.1cm]
    \text{where } \xx =& \vert\xx\vert(\cos(\psi(\xx)),\sin(\psi(\xx))),\nonumber\\[0.1cm]
    \phi_0(\xx) =& 4.2398+0.816735\cos(\psi(\xx)-0.2)
    -1.24397865\sin(2\psi(\xx)+0.1), \nonumber\\[0.1cm]
    \phi_1(\xx) =& 0.78167\sin(\psi(\xx)+0.5)- 2.24397865\cos(3\psi(\xx)-0.3) \nonumber\\[0.1cm]
    \phi_2(\xx) =& 1.127+1.2134875\cos(\psi(\xx)-0.65) -1.24397865\sin(2\psi(\xx)+0.1), \nonumber\\[0.1cm]
    \phi_3(\xx) =& 0.77-1.29\cos(4\psi(\xx)-0.35)+0.987\sin(2\psi(\xx)+0.14), \nonumber\\[0.1cm]
    r(\xx) =&  1.2927-0.929\cos(\psi(\xx)+0.34)+0.712\sin(3\psi(\xx)+0.14)\nonumber\\ &+\log(\vert\xx\vert+1.3),\nonumber\\[0.1cm]
    v(\xx) =& \left(1.1+\Re\left(H_{\vert\xx\vert^2 +1}^{(1)}(3)\right)\right)\exp\left(-\vert\xx-(0.027,\,0.0197)\vert^8\right)\nonumber\\ &\cdot(0.5+\sin(\xx_1(\xx_2-1))).\nonumber
\end{align}
In Figure~\ref{fig:conv2D-s-gen} we plot the difference {between values
obtained with grid sizes $h$ and $h/1.5$}
for the four different quadratures $\mathcal{U}_h^p$ \eqref{eq:Qp-general-explicit} and the punctured trapezoidal rule $T_h^0$. The order of accuracy shown is $1$ for the punctured trapezoidal rule and $p$ for the quadrature $\mathcal{U}_h^p$, which is what was expected. The error constant is determined by the value of $(\alpha,\beta)$. In all our tests we fixed $(\alpha,\beta)=(0.81,0.46)$. The stencils used for the different quadratures $Q_h^{k}$ needed to compose $\mathcal{U}_h^p$ are the same as the previous test, represented in Figure~\ref{fig:all_correction_grid_plot}. 

\subsection{Evaluating the layer potentials in the IBIM formulation}
\label{sub:num_3D}
We demonstrate the convergence and accuracy of the proposed quadrature rules by evaluating the single-layer, double-layer, and double-layer conjugate potentials with some smooth density $\rho$ on the surface $\Gamma\subset\mathbb{R}^3$:
\begin{equation*}
    \int_\Gamma G_0(\bar\xx^*,\bar{\mathbf{y}})\rho(\bar{\mathbf{y}})\dd\sigma_{\bar{\mathbf{y}}}  ,\  \int_\Gamma \dfrac{\partial G_0}{\partial\n_y}(\bar\xx^*,\bar{\mathbf{y}})\rho(\bar{\mathbf{y}})\dd\sigma_{\bar{\mathbf{y}}}  ,\  \int_\Gamma \dfrac{\partial G_0}{\partial\n_x}(\bar\xx^*,\bar{\mathbf{y}})\rho(\bar{\mathbf{y}})\dd\sigma_{\bar{\mathbf{y}}},\  \bar\xx^*\in\Gamma.
\end{equation*}

The integrals are first extended to the tubular neighborhood of  $T_\varepsilon$, as in \eqref{eq:general_S} using the compactly supported $C^{\infty}$ averaging function 
\begin{equation}\label{eq:averagingfunction}
    \delta(\eta)=
    \begin{cases}
        a\,\exp\left( \dfrac{2}{\eta^2-1} \right), & \text{ if } \vert\eta\vert<1, \\
        0, & \text{ otherwise};
    \end{cases}
\end{equation}
here $a\approx {{7.51393}}$ normalizes the integral $\int_{\mathbb{R}} \delta(\eta) \text{d}\eta$ to 1.

\subsubsection*{A numerical study on a smooth surface}
The surface chosen for the tests is a torus, centered {at a} randomly chosen point in 3D and rotated with randomly chosen angles along the $x$-, $y$- and $z$-axes. This is to avoid any symmetry of the uniform Cartesian grid which can influence the convergence behavior.
{This setup includes all the essential difficulties one may encounter when applying the proposed method to a smooth surface: non-convexity, finite reach from the geometry, and asymmetry in the discretized system. }

The torus is described by the following parametrization
\begin{equation}\label{eq:torus}
    \mathcal{T}(\theta,\phi)=Q\left(
    \begin{matrix}
    (R_2\cos\theta+R_1)\cos\phi \\
    (R_2\cos\theta+R_1)\sin\phi \\
    R_2\sin\theta
    \end{matrix}
    \right)+\mathbf{C}
\end{equation}
where $R_1=0.7$, $R_2=0.2$, $\mathbf{C}$ imposes a translation, and $Q=Q_z(c)Q_y(b)Q_x(a)$ is the composition of three rotation matrices; $Q_x(a)$, $Q_y(a)$, and $Q_z(a)$ are the matrices corresponding to a rotation by an angle $a$ around the $x$, $y$, and $z$ axes respectively. The parameters used for the translation and the rotations were: 
\begin{align*}
\mathbf{C}&=\big(0.5475547095598521,\, 0.6864792402110276,\, 0.3502726366462485\big)\cdot 10^{-1},\\
a&={0.199487}\cdot {10}^{{1}}, \\ b&={0.2540979476510170}\cdot {10}^{{1}}, \\ c&={0.4219760487439292}\cdot {10}^{{1}}.
\end{align*}
The known density function $\rho$ used in the test is defined using the parametrization of the torus:
\[
\rho(\bar{\mathbf{y}}) = \rho(\theta,\phi) = 1.38 + 2.196\sin\theta - 0.29837\cos\phi\,\sin\theta + 1.128\sin\phi\,\cos\theta\,.
\]
We present the errors
\begin{equation}\label{eq:error-Eh-SL-DL-DLC}
    \begin{array}{rl}
        E_{SL}^3(h) &=\, \left\vert \mathcal{V}^3_h\left[ G_0(\bar\xx^*,\bar{\mathbf{y}})\rho(\bar{\mathbf{y}}) \right] - \mathcal{V}^3_{h_{\min}}\left[ G_0(\bar\xx^*,\bar{\mathbf{y}})\rho(\bar{\mathbf{y}}) \right]  \right\vert, \\[0.3cm]
        E_{DL}^3(h) &=\, \left\vert \mathcal{V}^3_h\left[ \dfrac{\partial G_0}{\partial \n_y}(\bar\xx^*,\bar{\mathbf{y}})\rho(\bar{\mathbf{y}}) \right] - \mathcal{V}^3_{h_{\min}}\left[ \dfrac{\partial G_0}{\partial \n_y}(\bar\xx^*,\bar{\mathbf{y}})\rho(\bar{\mathbf{y}}) \right]  \right\vert, \\[0.3cm]
        E_{DLC}^3(h) &=\, \left\vert \mathcal{V}^3_h\left[ \dfrac{\partial G_0}{\partial \n_x}(\bar\xx^*,\bar{\mathbf{y}})\rho(\bar{\mathbf{y}}) \right] - \mathcal{V}^3_{h_{\min}}\left[ \dfrac{\partial G_0}{\partial \n_x}(\bar\xx^*,\bar{\mathbf{y}})\rho(\bar{\mathbf{y}}) \right]  \right\vert,
    \end{array}
\end{equation}
computed for a sequence of grid size values $\{h_i\}_i$,
where we used as reference value half of the smallest grid size $h_{\min}=\frac{1}{2}\min_i h_i$. We tested our third order rule $\mathcal{V}_h^3$ \eqref{eq:Q-3D-general}. Moreover we compared with the previously developed second order rule, denoted by $\mathcal{V}_h^2$, from~\cite{izzo2021corrected}.

In the presented simulations, we take the component $n_z$ of $\nn$ to be dominant if $\vert\tan\theta\vert<\sqrt{2}$, where $\nn/\vert\nn\vert=(\sin\theta\cos\phi,\sin\theta\sin\phi,\cos\theta)$. If instead $\vert\tan\theta\vert\ge\sqrt{2}$ and  $\vert\tan\phi\vert\ge 1$, we take $n_y$ to be dominant, and if $\vert\tan\theta\vert\ge\sqrt{2}$ and  $\vert\tan\phi\vert<1$ we take $n_x$ to be dominant. We used $\theta$ and $\phi$ to determine the dominant direction because of their extensive use in the rest of the code.

At each target point $\bar\xx^*$, the total error is the sum of the errors of the two-dimensional rule applied on each plane.
Recall that under the IBIM formulation, the kernel is singular along the surface's normal line passing through $\bar\xx^*$, 
and the singularity of the kernel on each plane lies at the intersection of the surface normal line and that plane.  
Since the normal lines of the surface generally do not align with the grid, the position of the singular point relative to the grid tends not to lie on any grid node. 
Recall further that the parameters $\alpha,\beta$ are used to described the position of the singular point relative to the closest grid node on the plane, and the 
error constants depend on them.  
Those parameters may change abruptly between planes, depending on which grid node in the plane 
is closest to the singular point. 
The closest grid nodes to each surface normal line certainly are expected to exhibit jumps as one refines the grids (decreases $h$).
Thus, as noted in \cite{izzo2021corrected}, the errors \eqref{eq:error-Eh-SL-DL-DLC} as functions of $h$ are generally  not smooth. Consequently we cannot see a clear slope.

\begin{figure}[h!]
    \begin{center}
    \includegraphics[scale=0.9]{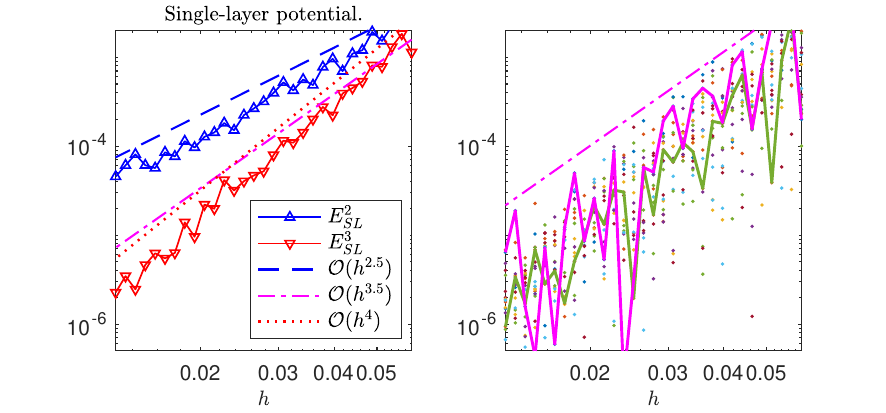}\\
    \includegraphics[scale=0.9]{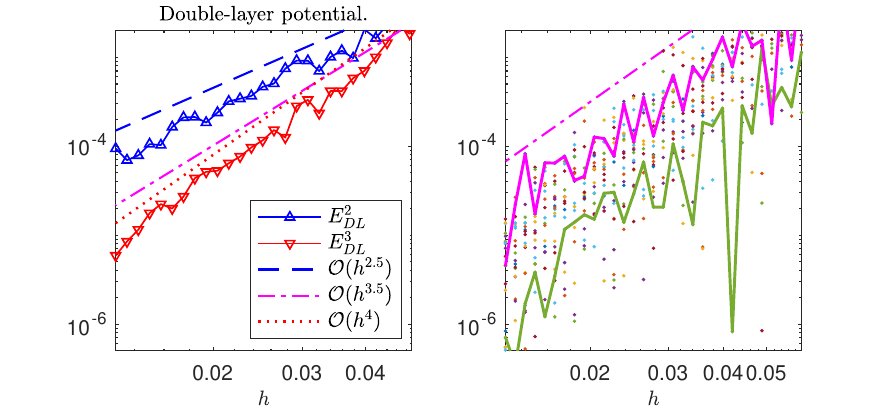}\\
    \includegraphics[scale=0.9]{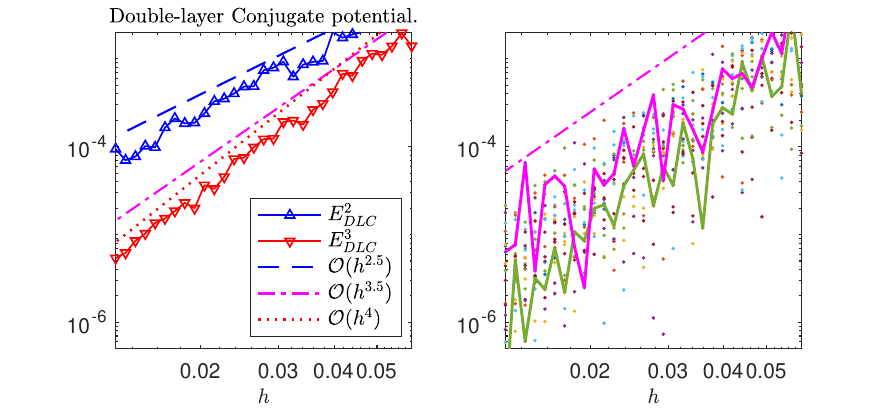}
    \end{center}
    \caption{\textbf{Errors in the evaluation of the three Laplace layer potentials}. The errors \eqref{eq:error-Eh-SL-DL-DLC} are computed for 20 randomly chosen target points on a tilted torus. The plots in the left column show the mean of the 20 errors. The plots in the right column show the scatter plot of the 20 target points. In the right plots we additionally highlight the behavior of two specific target points, to showcase a ``bad'' error (magenta line) and an ``average'' error (green line).}
    \label{fig:convergence3D}
\end{figure}

To show the overall convergence behavior we average the errors, defined in \eqref{eq:error-Eh-SL-DL-DLC},  over 20 target points, randomly chosen. The results can be seen in Figure~\ref{fig:convergence3D}. In the left column we present the averaged errors. 
In the right column we present a scatter plot of the errors at all the target points. We additionally highlight the errors corresponding to two specific target points to showcase an ``average'' error behavior (green line) and a ``bad'' error behavior (magenta line).

By construction of the quadrature rule \eqref{eq:Q-3D-general} we expect it to be third order accurate in $h$. However from the plots we observe order of accuracy $\ge 3.5$. 
We conjecture that an additional cancellation of errors occurs when adding the results from each plane (see \cite{izzo2021corrected} for a related discussion regarding $\mathcal{V}^2_h$). A rigorous analysis of this behavior is beyond the scope of this article.

\begin{figure}
    \begin{center}
    \includegraphics[scale=1]{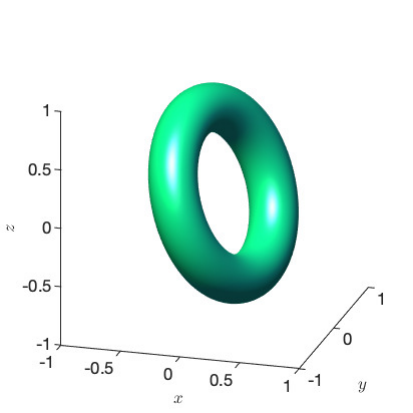}\includegraphics[scale=1]{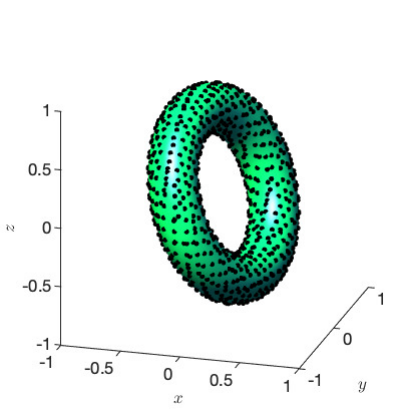}
    \end{center}
    \caption{\textbf{Torus test surface}. Left: the torus used in the tests. Right: the torus and the projections of the Cartesian grid nodes inside the tubular neighborhood $T_\varepsilon$. The projected nodes serve as the quadrature nodes.}
    \label{fig:tiltedtorusplot}
\end{figure}

Of course to test our algorithms, we retain no information about the parametrizations. The test torus is represented only by $d_\Gamma$ and $P_\Gamma$ on the given grid.
Figure~\ref{fig:tiltedtorusplot} shows the torus that we use and the points used in the quadrature rule for a given grid configuration. 
We use fourth-order centered differencing of $P_\Gamma$ on the grid to approximate the Jacobian $J_\Gamma$ (see \cite{tsaikublik16}). We also use fourth-order centered differencing of $P_\Gamma$ to find the third derivatives of $f$ needed for the functions $B$ and $C$ in \eqref{eq:surface-B-C}, as they are related via a linear system (see Appendix~\ref{sec:appendixB}).

%\vspace{0.5cm}
\subsubsection*{A numerical study on a more complicated surface}

\begin{figure}
    \begin{center}
        \includegraphics[scale=1.2]{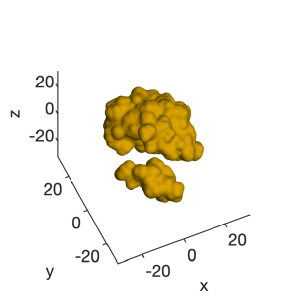}
    \end{center}
    \caption{{ \textbf{A solvent-molecule interface}. The surface is computed by the VISM method for biomolecular system p53-MDM2 (PDB ID 1YCR) \cite{kussie1996structure} from the Protein Data Bank (PDB). }}
    \label{fig:proteins}
\end{figure}

We present a test of the quadrature rule applied to IBIM for a more complicated surface, shown in Figure~\ref{fig:proteins}.
The surface represents the solvent-molecule interface of a complex biomolecular system immersed in a solvent \cite{zhang2021coupling}. A level set representation of the surface is generated using VISM \cite{wang2012level} by the authors of \cite{zhang2021coupling} on a $512^3$ Cartesian grid. 
We compute the relative error in the double-layer identity 
\[
\int_\Gamma \frac{\partial G_0}{\partial \n_y}(\bar{\xx},\bar{\yy})\dd\sigma_{\bar{\yy}} = - \frac12, \quad \bar{\xx}\in\Gamma,
\]
using the proposed method.
The relative error is defined as:
\begin{equation*}\label{eq:error-DL-identity-complex-geometry}
\mathcal{E}^p(h) = 2\sqrt{ h^3 \sum_{{\xx}_j\in h\mathbb{Z}^3 \cap T_\epsilon}  \left\vert \mathcal{V}_h^p\left[ \frac{\partial G_0}{\partial \n_y}(\bar{\xx}_j,\bar{\yy}) \right]  + \frac12 \right\vert^2\delta_{\Gamma,\epsilon}(\xx_j)},~~~\varepsilon=2h,
\end{equation*}
with $p=0$ denoting the punctured trapezoidal rule.

{ In our setup, inherited from the shared data set, the signed distance function is accurate up to distance $ \approx 9h$, which is minimally adequate for the application of $\mathcal{V}^2_h$.
The computed values are $\mathcal{E}^0(h)= 0.0159$ for the punctured trapezoidal rule and $\mathcal{E}^2(h)= 0.00142$ for the second order corrected rule. Furthermore, when applying $\mathcal{V}^3_h$, we notice that the resulting pointwise errors oscillate across grid nodes, $\bar\xx_j$, and do not appear to be smaller than those computed by $\mathcal{V}^2_h$. On some $\bar\xx_{j^\prime}$, the error even appear to be larger that those computed by the punctured trapezoidal rule. This is \emph{expected} because the grid does not yet resolve the fine geometry in this surface (notice in particular the narrow separation of the two connected components). }

%\newpage
% \backmatter
\section*{Declarations}

% \bmhead{Funding}
\subsection*{Funding}

Tsai's research is supported partially by National Science Foundation Grants DMS-1913209 and DMS-2110895.

% \bmhead{Conflict of interest}
\subsection*{Conflict of interest}

The authors declare that they have no known competing financial interests or personal relationships that could have appeared to influence the work reported in this paper.

% \bmhead{Authors' contribution}
\subsection*{CRediT taxonomy of authors' contribution}

Conceptualization: Olof Runborg, Richard Tsai; 
Methodology: Federico Izzo;
Investigation: Federico Izzo;
Software: Federico Izzo; 
Visualization: Federico Izzo;
Writing - original draft preparation: Federico Izzo, Olof Runborg, Richard Tsai;
Writing - review and editing: Federico Izzo, Olof Runborg, Richard Tsai;
Funding acquisition: Richard Tsai;
Supervision: Olof Runborg, Richard Tsai.

\printbibliography

\appendix
\section{Proofs of the lemmas and theorems}
\label{sec:appendix}

In this Section, we will prove the Lemmas and Theorems mentioned in Section~\ref{sec:corrected-trapezoidal-rules} and Section~\ref{sec:application-layer-potentials-3D}. 

\subsection{Proof of Theorem~\ref{thm:punctured-tr-s-ell-around-singularity}}
\label{sec:A:thm-s-singularity}

Consider a cut-off function $\psi\in C^\infty_c({\mathbb R}^n)$ such that
\begin{equation}\label{eq:psi-function}
\psi(\xx)=
\begin{cases}
1\,, & \vert\xx\vert\leq \frac{1}{2}\,,\\
0\,, & \vert\xx\vert\geq 1\,.
\end{cases}
\end{equation}
Then we can write $f$ as
\begin{align*}
    f(\xx) =& s(\xx)v(\xx) = s(\xx)v(\xx)\psi(\xx/r_0)+s(\xx)v(\xx)(1-\psi(\xx/r_0))\\ 
    =& \vert\xx\vert^j\ell(\vert\xx\vert,\xx/\vert\xx\vert)\psi(\xx/r_0)v(\xx)+s(\xx)v(\xx)(1-\psi(\xx/r_0))\\
    =& \vert\xx\vert^j\ell_1(\vert\xx\vert,\xx/\vert\xx\vert)v(\xx)+s(\xx)v(\xx)(1-\psi(\xx/r_0)).
\end{align*}
The first term is a function compactly supported in $B_{r_0}$, so by extending it to zero in ${\mathbb{R}}^n$ it satisfies the hypotheses of Theorem~\ref{thm:A:puncturederr}. Hence the result is valid for the first term.

The second term has regularity $C^\infty_c({\mathbb{R}}^n)$ and is zero in $B_{r_0/2}$, so the error for the punctured trapezoidal rule will decrease faster than any polynomial of $h$.

By combining the results for the two terms, we prove the result.
\subsubsection{Results on which Theorem~\ref{thm:punctured-tr-s-ell-around-singularity} depends.}
\label{sec:A:thm-s-singularity-supplementary}
\begin{theorem}\label{thm:A:puncturederr}
Suppose $v\in C^\infty_c({\mathbb R}^n)$ and $\ell\in C^{\infty}({\mathbb{R}}\times\s^{n-1})$.
Then, for integers $j\geq 1-n$, 
\begin{equation}\label{eq:A:theorem-statement-TR0-order}
    \left\vert\int_{{\mathbb{R}}^n} s(\xx)v(\xx) \text{\emph{d}}\xx - T^0_{h,\,\N_h}[s\, v]\right\vert \leq C h^{j+n}\,  ,\ \ 
  s(\xx)=\vert\xx\vert^{j}\, \ell\left(\vert\xx\vert,\frac{\xx}{\vert\xx\vert}\right)\,,
\end{equation}
where the constant $C$ is independent of $h$, but depends on $j$, $\ell$ and $v$.
\end{theorem}
\begin{proof}
Define $f(\xx):=\vert\xx\vert^j\ell(\vert\xx\vert,\xx/\vert\xx\vert)v(\xx)$, and consider the cut-off function $\psi\in C^\infty_c({\mathbb R}^n)$ \eqref{eq:psi-function}. 
Then we can write the punctured trapezoidal rule as 
\[
T^0_{h,\,\N_h}[f]=T_h[f(\,\cdot\,)(1-\psi(\cdot/h))],
\]
where we cut out the singularity point by multiplying by $1-\psi$ around $\zero$; the scaling by $h$ ensures that, for fixed $h$, only the node in the singularity point is cut out. 
This allows us to split the
error of the punctured trapezoidal rule 
as
\begin{align*}
\int_{{\mathbb{R}}^n}f(\xx)\dd\xx-T_{h,\,\N_h}^0[f] =&\, \underbrace{\int_{{\mathbb{R}}^n}f(\xx)\psi(\xx/h)\dd\xx}_{(\textbf{I})} \\ &+ \underbrace{\int_{{\mathbb{R}}^n}f(\xx)(1-\psi(\xx/h))\dd\xx-T_h[f(\cdot)(1-\psi(\cdot/h))]}_{(\textbf{II})}\,.
\end{align*}
We will consider the two terms (\textbf{I}), (\textbf{II}) separately, and prove that both can be bounded by $Ch^{j+n}$.\\

\noindent(\textbf{I}): Given the compact support of $\psi$, the 
integral is reduced to an integral
over
$\{\vert\xx\vert\leq h\}$:
\begin{align*}
    \int_{{\mathbb{R}}^n}f(\xx)\psi(\xx/h)\dd\xx =&\, \int_{\vert\xx\vert\leq h}v(\xx)\vert\xx\vert^j\ell(\vert\xx\vert,\xx/\vert\xx\vert)\psi(\xx/h)\dd\xx \\
    =&\, h^{j+n}\int_{\vert\xx\vert\leq 1}v(h\xx)\vert\xx\vert^j\ell(\vert h\xx\vert,\xx/\vert\xx\vert)\psi(\xx)\dd\xx \\
    \Rightarrow\ \left\vert\int_{{\mathbb{R}}^n}f(\xx)\psi(\xx/h)\dd\xx\right\vert \leq&\, h^{j+n} \vert v\vert_\infty\vert\ell\vert_\infty \int_{\vert\xx\vert\leq 1}\vert\xx\vert^j\dd\xx \leq C_1 h^{j+n},
\end{align*}
since $\vert\xx\vert^j$ is integrable
as $j\geq 1-n$.
We have proven the estimate
for the first term.\\

\noindent(\textbf{II}): For the second term, knowing that the volume of the fundamental parallelepiped of the lattice $V:=(h\mathbb{Z})^n$ is $h^n$ and that the dual lattice is $V^*=(h^{-1}\mathbb{Z})^n$, we use the Poisson summation formula:
\[
T_h[f] = h^n\sum_{\mathbf{j}\in V} f(\mathbf{j}) = \dfrac{h^n}{h^n} \sum_{\mathbf{l}\in V^*} \hat f\left( {\mathbf{l}} \right) = \int_{{\mathbb{R}}^n} f(\xx)\dd\xx +  \sum_{\kk\neq\zero} \hat f\left( \dfrac{\kk}{h} \right)\,.
\]
Then the error in (\textbf{II}) is:
\[
T_h\left[ f(\cdot)\,(1-\psi(\cdot/h)) \right] - \int_{{\mathbb{R}}^n}f(\xx)(1-\psi(\xx/h))\dd\xx = \sum_{\kk\neq\zero} \hat f_\psi (\kk,h),
\]
where
\begin{align*}
\hat f_\psi(\kk,h) :=& \hat f(\kk/h) = \int_{{\mathbb{R}}^n} f(\xx)(1-\psi(\xx/h))e^{-2\pi\ii\kk \cdot \xx/h}\dd\xx \\ =& h^n\int_{{\mathbb{R}}^n} f(h\xx)(1-\psi(\xx))e^{-2\pi\ii\kk\cdot\xx}\dd\xx\,.
\end{align*}
Using integration by parts separately on each of the variables, we find
\begin{align*}
    \int_{{\mathbb{R}}^n} \partial_{\xx}^{\beta}[f(h\xx)(1-\psi(\xx))]e^{-2\pi\ii\kk\cdot\xx}\dd\xx &= 2\pi\ii\,k_j \int_{{\mathbb{R}}^n} \partial_\xx^{\beta-e_j}[f(h\xx)(1-\psi(\xx))]e^{-2\pi\ii \kk\cdot\xx}\dd\xx \\
    &= (2\pi\ii\kk)^{\beta} \int_{{\mathbb{R}}^2}f(h\xx)(1-\psi(\xx))e^{-2\pi\ii\kk\cdot\xx}\dd\xx.
\end{align*}
For the Laplacian operator applied $q$ times we therefore have
\begin{align*}
    & \int_{{\mathbb{R}}^n} \Delta^{q}[f(h\xx)(1-\psi(\xx))]e^{-2\pi\ii\kk\cdot\xx}\dd\xx \\
    =& -4\pi^2\left(\sum_{j=1}^n k_j^2\right)\int_{{\mathbb{R}}^n} \Delta^{q-1}[f(h\xx)(1-\psi(\xx))]e^{-2\pi\ii\kk\cdot\xx}\dd\xx \\
    =& (-1)^q (2\pi)^{2q} \vert\kk\vert^{2q} \int_{{\mathbb{R}}^n}f(h\xx)(1-\psi(\xx))e^{-2\pi\ii\kk\cdot\xx}\dd\xx \,.
\end{align*}
We use this result to find an expression we can bound using Lemma \ref{lem:A:hp-hbeta-2-ell(r,s)}; given an integer $q$, we find
\begin{align*}
    \left\vert \hat f_\psi(\kk,h) \right\vert &\leq \dfrac{h^n}{(2\pi)^{2q}\vert\kk\vert^{2q}} \int_{{\mathbb{R}}^n}\left\vert \Delta^q[f(h\xx)(1-\psi(\xx)]e^{\ii\,\kk\cdot\xx} \right\vert \dd\xx  \\
    &\leq \dfrac{h^n}{(2\pi)^{2q}\vert\kk\vert^{2q}} \sum_{\vert\beta\vert=2q} c_\beta \int_{{\mathbb{R}}^n} \left\vert \partial_\xx^\beta[f(h\xx)(1-\psi(\xx)] \right\vert \dd\xx  \\
    &\leq \dfrac{h^n}{(2\pi)^{2q}\vert\kk\vert^{2q}} \sum_{\vert\beta\vert=2q} \tilde c_\beta (h^j+h^{\vert\beta\vert-n}) = \bar c_\beta \dfrac{h^{j+n}+h^{2q}}{\vert\kk\vert^{2q}}\,.
\end{align*}
Then the series of Fourier coefficients is
\begin{align*}
   \left\vert T_h\left[ f(\cdot)\,(1-\psi(\cdot/h)) \right] - \int_{{\mathbb{R}}^n}f(\xx)(1-\psi(\xx/h))\dd\xx \right\vert &\leq  \sum_{\kk\neq\zero} \left\vert \hat f_\psi(\kk,h) \right\vert  \\
   & \leq \bar c_\beta\sum_{\kk\neq\zero} \dfrac{h^{j+n}+h^{2q}}{\vert\kk\vert^{2q}}\,.
\end{align*}
The series converges if $2q>n$, and the leading order is $h^{j+n}$ if $2q\geq j+n$, so by taking $q\geq \max(1+n/2,(n+j)/2)$, we find the result sought. Combining the results for (\textbf{I}) and (\textbf{II}), we find the bound
\[
\left\vert \int_{{\mathbb{R}}^n}f(\xx)\dd\xx-T_{h,\,\N_h}^0[f] \right\vert \leq C_\beta h^{j+n}\,. 
\]
This proves the theorem.
\end{proof}

We use the notation $\xx=(x_1,x_2,\dots,x_n)=\sum_{l=1}^n x_l e_l$, and indicate with $e_l$ the $l$-th element of the standard ${\mathbb{R}}^n$ basis.
\begin{lemma}\label{lem:A:hp-hbeta-2-ell(r,s)}
Let $g,\psi\in C^\infty_c({\mathbb{R}}^n)$, $\ell\in C^\infty({\mathbb{R}}\times\s^{n-1})$, where $\psi$ is such that 
\[
\psi(\xx)=
\begin{cases}
1 & \vert\xx\vert\leq \frac{1}{2}\,,\\[0.1cm]
0 & \vert\xx\vert\ge 1\,.
\end{cases}
\]
Let $j\geq 1-n$, and $f(\xx)=\vert\xx\vert^j \ell(\vert\xx\vert,\xx/\vert\xx\vert)g(\xx)$; then, for any multi-index $\beta\in\mathbb{N}^n_0$ it exists a constant $C_\beta$ independent of $h$ such that, for $0<h\leq 1$,
\begin{equation}\label{eq:lemma-thm-thesis}
    \int_{{\mathbb{R}}^2}\Big\vert \partial_\xx^\beta\left[ f(h\xx)(1-\psi(\xx)) \right] \Big\vert\text{\emph{d}} \xx\leq C_\beta(h^{j}+h^{\vert\beta\vert-n})\,.
\end{equation}
\end{lemma}
\begin{proof}
Given $\beta\in\mathbb{N}^n_0$, we first prove that there exist functions $f_{\beta}:{\mathbb{R}}\times\s^{n-1}\to{\mathbb{R}}$
in $C_c^\infty({\mathbb{R}}\times\s^{n-1})$ such that
\begin{equation}\label{eq:lemma-thm-firststep}
    \partial_\xx^\beta f(\xx)=\vert\xx\vert^{j-\vert\beta\vert} f_{\beta}(\vert\xx\vert,\xx/\vert\xx\vert)\,.
\end{equation}
We prove this by induction. The induction base $\beta=\zero$ is true because
\[
\partial_\xx^\zero f(\xx)=f(\xx)=\vert\xx\vert^j \ell(\vert\xx\vert,\xx/\vert\xx\vert)g(\xx)=:\vert\xx\vert^j f_{\zero}(\vert\xx\vert,\xx/\vert\xx\vert),
\]
where $f_{\zero}\in C^\infty_c({\mathbb{R}}\times\s^{n-1})$. For the induction step we assume that \eqref{eq:lemma-thm-firststep} is true for $\beta$ and prove it for $\beta+e_l$:
\[
\partial_\xx^{\beta+e_l}f(\xx)=\partial_\xx^{e_l} \vert\xx\vert^{j-\vert\beta\vert} f_{\beta}(\vert\xx\vert,\xx/\vert\xx\vert)\,.
\]
By computing the derivative we find
\begin{align*}
    \partial_\xx^{e_l}\vert\xx\vert^{j-\vert\beta\vert} f_{\beta}\left(\vert\xx\vert,\dfrac{\xx}{\vert\xx\vert}\right) =&\, \vert\xx\vert^{j-\vert\beta\vert-1}\Bigg[(j-\vert\beta\vert)\,\left(\dfrac{\xx}{\vert\xx\vert}\right)_l f_{\beta}\left(\vert\xx\vert,\dfrac{\xx}{\vert\xx\vert}\right) \\
    &+ \nabla_{\mathbf{u}}f_{\beta}\left(\vert\xx\vert,\dfrac{\xx}{\vert\xx\vert}\right)\cdot\left(e_l-\left(\dfrac{\xx}{\vert\xx\vert}\right)_l\dfrac{\xx}{\vert\xx\vert}\right) \\
    &+ \vert\xx\vert\left(\dfrac{\xx}{\vert\xx\vert}\right)_l\partial_r f_{\beta}\left(\vert\xx\vert,\dfrac{\xx}{\vert\xx\vert}\right)\Bigg]\\
    =&:\vert\xx\vert^{j-\vert\beta\vert-1}f_{\beta+e_l}\left(\vert\xx\vert,\dfrac{\xx}{\vert\xx\vert}\right)\,.
\end{align*}
Because $f_\beta\in C^\infty_c({\mathbb{R}}\times\s^{n-1})$ the same is also true for $f_{\beta+e_l}$.\\
The next step is to expand the derivative in \eqref{eq:lemma-thm-thesis} and use \eqref{eq:lemma-thm-firststep}, and then bound it:
\begin{align*}
    \partial_\xx^\beta\left[ f(h\xx)(1-\psi(\xx)) \right] =&
     \sum_{\nu\leq\beta}{\beta\choose\nu}\partial^{\beta-\nu}[1-\psi(\xx)]\,h^{\vert\nu\vert}\partial^\nu f(h\xx) \\
     =& \sum_{\nu\leq\beta}{\beta\choose\nu}\partial^{\beta-\nu}[1-\psi(\xx)]h^{j}\vert\xx\vert^{j-\vert\nu\vert} f_{\nu}(\vert h\xx\vert,\mathbf{u})\,.
\end{align*}
We use the properties of $\psi$, and the compact support of $f_\nu$. Let $L>0$ be such that $\forall \nu\leq\beta$, supp$\,f_\nu$ is contained in the ball $B_L(\zero)$. Note furthermore that the derivatives of $\psi$ are compactly supported in the annulus $\{\xx\in{\mathbb{R}}^n\,:\,\frac{1}{2}\leq\vert\xx\vert\leq 1\}$. From this we can say that
\[
\Big\vert \partial_\xx^\beta\left[ f(h\xx)(1-\psi(\xx)) \right] \Big\vert\leq\ C
\begin{cases}
0\,, & \vert\xx\vert\leq \frac{1}{2}\,, \\
h^{j}\,, & \frac{1}{2}\leq\vert\xx\vert\leq 1\,, \\
h^{j}\vert\xx\vert^{j-\vert\beta\vert}\,, & 1\leq\vert\xx\vert\leq L/h\,, \\
0\,, & \vert\xx\vert>L/h\,.
\end{cases}
\]
We use these bounds in the evaluation of the integral, and after passing to polar coordinates we arrive at \eqref{eq:lemma-thm-thesis} via
\begin{align*}
    \int_{{\mathbb{R}}^n}\Big\vert \partial_\xx^\beta\left[ f(h\xx)(1-\psi(\xx)) \right] \Big\vert\dd\xx \leq& C_1\int_{1/2}^1 h^j r^{n-1}\dd r+C_2 h^j \int_{1}^{L/h}r^{j-\vert\beta\vert+n-1}\dd r \\
    =& \bar C_1 h^j + C_2 h^{\vert\beta\vert-n} \int_{h}^L r^{j-\vert\beta\vert+n-1}\dd r \\
    =& \bar C_1 h^j + C_2 h^{\vert\beta\vert-n} \left( C_3+C_4 h^{j-\vert\beta\vert+n} \right) \\
    \leq & C_\beta\,( h^j +  h^{\vert\beta\vert-n} )\,.
\end{align*}
The lemma is proven.
\end{proof}

\subsection{Proof of Lemma~\ref{lem:f-over-distance-expansion}}
\label{sub:A:f-over-distance-expansion}
For any $\mathbf{u}\in\s^1$, we expand $\ell$ around $r=0$ and write the remainder in integral form:
\[
\ell(r,\mathbf{u})=\sum_{j=0}^q\dfrac{1}{j!}\partial_r^j\ell(0,\mathbf{u})r^j+\dfrac{r^{q+1}}{q!}\int_0^1\partial_r^{q+1}\ell(tr,\mathbf{u})(1-t)^q\dd t\,.
\]
Then 
\begin{align*}
\triangle_q s(\xx)&=\dfrac{1}{\vert\xx\vert}\ell\left(\vert\xx\vert,\dfrac{\xx}{\vert\xx\vert}\right)-\sum_{j=0}^q \dfrac{1}{j!}\partial_r^j \ell\left(0,\dfrac{\xx}{\vert\xx\vert}\right)\vert\xx\vert^{j-1}\\
&= \dfrac{\vert\xx\vert^q}{q!}\int_0^1(1-t)^q\partial_r^{q+1}\ell(t\vert\xx\vert,\xx/\vert\xx\vert)\dd t=\vert\xx\vert^q\sigma (\vert\xx\vert,\xx/\vert\xx\vert)\,,
\end{align*}
where $\sigma\in C^\infty((-r_0,r_0)\times \s^1)$ because $\ell\in C^\infty((-r_0,r_0)\times \s^1)$. The lemma is thus proven.

\subsection{Proof of Lemma~\ref{lem:properties-of-h}}
\label{sub:A:properties-of-h}
    The first two identities in \eqref{eq:hprop1} follows since
    $P_\Gamma(\bar\xx^*+(\zero,z')_B) = \bar\xx^*$ for all $z'$, as was already pointed out
    in Section~\ref{subsub:proj-mapping-to-layer-kernels}.
    For the second part, 
    we note that the surface normal at
    the point $\bar\xx^*+ \big( \yy_{\text{p}},f(\yy_{\text{p}})\big)_B$
    is parallell to $(-\nabla f(\yy_{\text{p}}),1)_B$.
    Therefore, there is a $t\in{\mathbb{R}}$
    such that
    $$
    \bar\xx^*+({\mathbf{y}}',z')_B=
    \bar\xx^*+ \big( \yy_{\text{p}},f(\yy_{\text{p}})\big)_B+t
    (-\nabla f(\yy_{\text{p}}),1)_B,
    $$
    which implies that
    \begin{equation}\label{eq:A:x-projection-general-f}
    {\mathbf{y}}'= \yy_{\text{p}}-(z'-f(\yy_{\text{p}})) \nabla f(\yy_{\text{p}})=:\F(\yy_{\text{p}}).
    \end{equation}
    Using the fact that $\yy_{\text{p}}=\bh({\mathbf{y}}',z')$
    and
    differentiating both sides with
    respect to ${\mathbf{y}}'$ gives us,
    \begin{align*}
    I&=
    \frac{\partial \F(\yp)}{\partial \yp}^T
    \frac{\partial \bh}{\partial {\mathbf{y}}}
    =\frac{\partial \bh}{\partial {\mathbf{y}}}
    -\left(
    (z'-f(\yy_{\text{p}}))
    \frac{\partial^2 f}{\partial {\mathbf{y}}^2}(\yy_{\text{p}})
    -\nabla f(\yy_{\text{p}})
    \nabla f(\yy_{\text{p}})^T\right)\frac{\partial \bh}{\partial {\mathbf{y}}},
    \end{align*}
    and the result follows upon evaluating
    at ${\mathbf{y}}'=\yy_{\text{p}}=\zero$ and using \eqref{eq:f-properties-def-M}.
    Since $\bh$ is smooth on ${\mathcal M}_L$
    the matrix $D$ must thus be well-defined.
    
    For the second order term in the Taylor expansion,
    we write $\yp=(y_1,y_2)$, $\bh=(h_1,h_2)^T$ and
    $\F=(F_1,F_2)^T$. We then
    get for $j=1,2$,
    $$
    \zero =\frac{\partial^2 F_j(\bh)}{\partial {\mathbf{y}}^2}
    =
    \frac{\partial F_j(\yp)}{\partial y_1}
    \frac{\partial^2 h_1}{\partial {\mathbf{y}}^2}+
    \frac{\partial F_j(\yp)}{\partial y_2}
    \frac{\partial^2 h_2}{\partial {\mathbf{y}}^2}
    +
    \frac{\partial\bh}{\partial {\mathbf{y}}}^T
    \frac{\partial^2 F_j(\yp)}{\partial \yp^2}
    \frac{\partial\bh}{\partial {\mathbf{y}}}.
    $$
    From the expressions above we have that
    $\frac{\partial \F(\zero)}{\partial \yp}=D^{-1}(z)$.
    Therefore, evaluating at ${\mathbf{y}}'=\yy_{\text{p}}=\zero$,
    yields
    $$
    \zero =
    D(z')^{-1}_{jj}
    \frac{\partial^2 h_j}{\partial {\mathbf{y}}^2}
    +
    D(z')^T
    \frac{\partial^2 F_j(\yp)}{\partial \yp^2}
    D(z'), \qquad j=1,2.
    $$
    Since
    $$
    \left.\frac{\partial^2 F_j(\yp)}{\partial \yp^2}
    \right\vert_{\yp=0}
    =\left.-z'\frac{\partial}{\partial y_j}\frac{\partial^2 f}{\partial \yp^2}\right\vert_{\yp=0},
    $$
    we finally get
    \begin{align*}
    \frac{1}{2}  \Vector{
      {{\mathbf{y}}'}^T\frac{\partial^2 h_1}{\partial {\mathbf{y}}^2}{\mathbf{y}}'\\
      {{\mathbf{y}}'}^T\frac{\partial^2 h_2}{\partial {\mathbf{y}}^2}{\mathbf{y}}'
      }
      =&
      -\frac{1}{2} D(z')\Vector{
     {{\mathbf{y}}'}^TD(z')^T
    \frac{\partial^2 F_1(\yp)}{\partial \yp^2}
    D(z'){{\mathbf{y}}'}\\
    {{\mathbf{y}}'}^TD(z')^T
    \frac{\partial^2 F_2(\yp)}{\partial \yp^2}
    D(z'){{\mathbf{y}}'}}
    \\
    =&\ z' D(z')C\Bigl(D(z'){\mathbf{y}}',D(z'){\mathbf{y}}'\Bigr).
    \end{align*}
    This gives \eqref{eq:hig_order_proj_x} and the lemma is proven.

\subsection{Proof of Lemma~\ref{lem:g=f-over-distance}}
\label{sub:A:g=f-over-distance}

For the first function, using the hypothesis $\bar\bg(\zero)=\zero$ and the notation $\xx=\vert\xx\vert\mathbf{u}$ with $\xx/\vert\xx\vert=:\mathbf{u}\in\s^{m-1}$ we write the expansion around $\xx=\zero$ as
\begin{align*}
\bar\bg(\xx) =& \bar\bg(\zero)+D\bar\bg(\zero)\xx+\sum_{\vert\nu\vert=2}E_{\bar\bg,\nu}(\xx)\xx^\nu \\ =& \vert\xx\vert\left( D\bar\bg(\zero)\mathbf{u}+\vert\xx\vert \sum_{\vert\nu\vert=2}E_{\bar\bg,\nu}(\xx)\mathbf{u}^\nu\right) =:\vert\xx\vert f(\vert\xx\vert,\mathbf{u})\,,
\end{align*}
where $E_{\bar\bg,\nu}(\xx):=\frac{2}{\nu!}\int_0^1(1-t)\partial^\nu \bar\bg(t\xx)\dd t$
is given by
the integral form of the remainder term.
Using the full rank of $D\bar\bg(\zero)$, there exists $0<r_1\leq r_0$ be such that $f(\vert\xx\vert,\mathbf{u})\neq 0$ in $(-r_1,r_1)\times \s^{m-1}$. Then 
\[
\dfrac{1}{\vert\bar\bg(\xx)\vert}=\dfrac{1}{\vert\xx\vert}\dfrac{1}{\vert f(\vert\xx\vert,\mathbf{u})\vert}=\dfrac{1}{\vert\xx\vert}\ell_1\left(\vert\xx\vert,\mathbf{u}\right)\,,
\]
and from the hypotheses on $D\bar\bg(\zero)$ and on the smoothness of $\bar\bg$, $\ell_1$ is $C^\infty((-r_1,r_1)\times \s^{m-1})$.

For the second function form, let $r(\xx):=\bar\bp(\xx)^T\bar\bg(\xx)$; then $\nabla r(\xx)=\bar\bg(\xx)^TD\bar\bp(\xx)+\bar\bp(\xx)^TD\bar\bg(\xx)$. Using the hypothesis $\bar\bp(\zero)^TD\bar\bg(\zero)=\zero$, we write the expansion of $r$ around $\xx=\zero$ using the integral form of the remainder:
\[
r(\xx)=r(\zero)+\nabla r(\zero)\xx+\sum_{\vert\nu\vert=2}E_{r,\nu}(\xx)\xx^\nu=\vert\xx\vert^2\sum_{\vert\nu\vert=2}E_{r,\nu}(\xx)\mathbf{u}^\nu\,,
\]
where $E_{r,\nu}(\xx):=\frac{2}{\nu!}\int_0^1(1-t)\partial^\nu r(t\xx)\dd t$, so that we find
\[
\dfrac{\bar\bp(\xx)^T\bar\bg(\xx)}{\vert\bar\bg(\xx)\vert^3}=\dfrac{\vert\xx\vert^2\sum_{\vert\nu\vert=2}E_{r,\nu}(\xx)\mathbf{u}^\nu}{\vert\xx\vert^3f(\vert\xx\vert,\mathbf{u})^3}=\dfrac{1}{\vert\xx\vert}\dfrac{\sum_{\vert\nu\vert=2}E_{r,\nu}(\xx)\mathbf{u}^\nu}{f(\vert\xx\vert,\mathbf{u})^3}=\dfrac{1}{\vert\xx\vert}\ell_2(\vert\xx\vert,\mathbf{u})\,.
\]
From the hypotheses on the smoothness of $\bar\bg$ and $\bar\bp$, $\ell_2$ is $C^\infty((-r_1,r_1)\times \s^{m-1})$ and the result is proven.

\section{Computation of the derivatives of the local surface function}
\label{sec:appendixB}
In this Section we will show how to find numerically the derivatives of $f$ in the Implicit Boundary Integral Methods setting of Section~\ref{sec:application-layer-potentials-3D}. The derivatives are needed to evaluate the functions $B$ and $C$ of~\eqref{eq:surface-B-C}, which are used in the approximated kernels \eqref{eq:s0s1-kernels}. 

The first derivatives and the mixed second derivatives are zero by construction, so we will show how to find the pure second derivatives and all the third derivatives.\\

Let $\bar{\mathbf{z}}$ be an arbitrary point in $T_\varepsilon$, and $\eta=d_\Gamma(\bar{\mathbf{z}})$. Let $\Gamma_\eta:=\{\bar{\mathbf{z}}\in T_\varepsilon\,:\,d_\Gamma(\bar{\mathbf{z}})=\eta\}$ be the surface parallel to $\Gamma$ at signed distance $\eta$. 

The pure second derivatives of $f$ at $P_\Gamma(\bar{\mathbf{z}})$, $f_{xx},f_{yy}$, are the principal directions $\ka_1,\ka_2$ of $\Gamma$ at $P_\Gamma(\bar{\mathbf{z}})$. We find the principal curvatures $g_1,g_2$ of $\Gamma_\eta$ in $\bar{\mathbf{z}}$ via the Hessian of $d_\Gamma$ at $\bar{\mathbf{z}}$:
\[
H_{{d}_\Gamma}(\bar{\mathbf{z}})=\nabla^{2}{d}_\Gamma(\bar{\mathbf{z}})=\left[\begin{array}{ccc}
\nn & \bar{\bm{\tau}}_{1} & \bar{\bm{\tau}}_{2}\end{array}\right]\begin{bmatrix}0\\
 & -g_{1}\\
 &  & -g_{2}
\end{bmatrix}\left[\begin{array}{ccc}
\nn & \bar{\bm{\tau}}_{1} & \bar{\bm{\tau}}_{2}\end{array}\right]^{T}
\]
where $\bar{\bm{\tau}}_1$, $\bar{\bm{\tau}}_2$ are the principal directions and $\nn$ is the normal to $\Gamma$ in $P_\Gamma(\mathbf{z})$. 
In practice, the values of either $P_\Gamma$ or $d_\Gamma$ are given on the grid nodes. The principal directions and curvatures are computed from eigendecomposition of third order numerical approximations of the Hessian, $H_{d_\Gamma}$. Alternatively, one can obtain this information from the derivative matrix of $P_\Gamma$, see \cite{tsaikublik16}.
Then the following relation lets us find the principal curvatures $\ka_1,\ka_2$ from $g_1,g_2$ and $\eta$:
\[
-\kappa_{i}=\frac{-g_{i}}{1+\eta g_{i}},~~~i=1,2.
\]

The third derivatives of $f$ can be found by computing the second derivatives with respect to ${\mathbf{y}}'$ of $\mathbf{h}({\mathbf{y}}',z')$ from Section~\ref{subsub:proj-mapping-to-layer-kernels}. By differentiating twice \eqref{eq:A:x-projection-general-f} with respect to ${\mathbf{y}}'=(x,y)$ with $\mathbf{h}({\mathbf{y}}',z')=\yy_{\text{p}}=(h_1,h_2)$ and evaluating in ${\mathbf{y}}'=\zero$, we find the following two linear systems:

\begin{align}
&V
\left(
\begin{array}{c}
f_{xxx} \\
f_{xxy} \\ 
f_{xyx} \\ 
f_{xyy}
\end{array}
\right) = \dfrac{1-z'\kappa_1}{z'}
\left(
\begin{matrix}
\frac{\partial^2 h_1}{\partial x^2}\\[0.1cm] 
\frac{\partial^2 h_1}{\partial x\partial y}\\[0.1cm] 
\frac{\partial^2 h_1}{\partial y\partial x}\\[0.1cm]
\frac{\partial^2 h_1}{\partial y^2}
\end{matrix}
\right), \label{eq:fxxx} 
\hspace{1cm}
V
\left(
\begin{array}{c}
f_{yxx} \\
f_{yxy} \\ 
f_{yyx} \\ 
f_{yyy}
\end{array}
\right) = \dfrac{1-z'\kappa_2}{z'}
\left(
\begin{matrix}
\frac{\partial^2 h_2}{\partial x^2}\\[0.1cm] 
\frac{\partial^2 h_2}{\partial x\partial y}\\[0.1cm] 
\frac{\partial^2 h_2}{\partial y\partial x}\\[0.1cm] 
\frac{\partial^2 h_2}{\partial y^2}
\end{matrix}
\right),
\end{align}
\begin{align*}
 \text{where }\ 
&V:=\left(
\begin{matrix}
\left(\frac{\partial h_1}{\partial x}\right)^2 & \frac{\partial h_1}{\partial x}\frac{\partial h_2}{\partial x} & \frac{\partial h_1}{\partial x}\frac{\partial h_2}{\partial x} & \left(\frac{\partial h_2}{\partial x}\right)^2\\
\frac{\partial h_1}{\partial x}\frac{\partial h_1}{\partial y} & \frac{\partial h_1}{\partial x}\frac{\partial h_2}{\partial y} & \frac{\partial h_1}{\partial y}\frac{\partial h_2}{\partial x} & \frac{\partial h_2}{\partial x}\frac{\partial h_2}{\partial y}\\
\frac{\partial h_1}{\partial x}\frac{\partial h_1}{\partial y} & \frac{\partial h_1}{\partial y}\frac{\partial h_2}{\partial x} & \frac{\partial h_1}{\partial x}\frac{\partial h_2}{\partial y} & \frac{\partial h_2}{\partial x}\frac{\partial h_2}{\partial y}\\
\left(\frac{\partial h_1}{\partial y}\right)^2 & \frac{\partial h_1}{\partial y}\frac{\partial h_2}{\partial y} & \frac{\partial h_1}{\partial y}\frac{\partial h_2}{\partial y} & \left(\frac{\partial h_2}{\partial y}\right)^2
\end{matrix}
\right).\nonumber
\end{align*}

We find the first and second derivatives of $\mathbf{h}({\mathbf{y}}',z')$ by computing the derivatives of $P_\Gamma$ in $\bar{\mathbf{z}}$ and applying a change of basis transformation. 

By construction $\bar{\mathbf{z}}=P_\Gamma(\bar{\mathbf{z}})+\eta \nn$. Then we use the closest point projections of the grid nodes around $\bar{\mathbf{z}}$, 
\[
\bar{\mathbf{v}}_{ijk}:=P_\Gamma(\bar{\mathbf{z}}+(i,j,k)h),\ \ i,j,k=-2,-1,0,1,2.
\]
In the $B$ basis, these points are expressed as $\bar{\mathbf{v}}_{ijk}=\bar\xx^*+(\bar{\mathbf{w}}_{ijk})_B$, where $\bar{\mathbf{w}}_{ijk}=Q^{-1}(\bar{\mathbf{v}}_{ijk}-\bar\xx^*)$. We apply finite differences (central differences of 4th order in this case) to the component of the nodes $\bar{\mathbf{w}}_{ijk}=(X_{ijk},Y_{ijk},Z_{ijk})$ to compute 
\[
W_1\approx\nabla X,\ \ W_2\approx\nabla Y,\ \ W_3\approx\nabla^2X,\ \ W_4\approx\nabla ^2Y.
\]
We can then use these approximations to find the derivatives of $h_i$, $i=1,2$ by applying the following transformations:
\begin{equation*}
\begin{array}{llll}
    \frac{\partial h_1}{\partial x} = \bar{\bm{\tau}}_1^T W_1, & \frac{\partial h_1}{\partial y} = \bar{\bm{\tau}}_2^T W_1 , & 
    \frac{\partial h_2}{\partial x} = \bar{\bm{\tau}}_1^T W_2, & \frac{\partial h_2}{\partial y} = \bar{\bm{\tau}}_2^T W_2, \\
    \frac{\partial^2 h_1}{\partial x^2} = \bar{\bm{\tau}}_1^TW_3\bar{\bm{\tau}}_1 , &  \frac{\partial^2 h_1}{\partial x\partial y} = \bar{\bm{\tau}}_2^TW_3\bar{\bm{\tau}}_1 , &  \frac{\partial^2 h_1}{\partial y\partial x} = \bar{\bm{\tau}}_1^TW_3\bar{\bm{\tau}}_2 , &  \frac{\partial^2 h_1}{\partial y^2} = \bar{\bm{\tau}}_2^TW_3\bar{\bm{\tau}}_2,\\
    \frac{\partial^2 h_2}{\partial x^2} = \bar{\bm{\tau}}_1^TW_4\bar{\bm{\tau}}_1 , &  \frac{\partial^2 h_2}{\partial x\partial y} = \bar{\bm{\tau}}_2^TW_4\bar{\bm{\tau}}_1 , &  \frac{\partial^2 h_2}{\partial y\partial x} = \bar{\bm{\tau}}_1^TW_4\bar{\bm{\tau}}_2 , &  \frac{\partial^2 h_2}{\partial y^2} = \bar{\bm{\tau}}_2^TW_4\bar{\bm{\tau}}_2.
\end{array}
\end{equation*} 

Finally, we solve the two systems \eqref{eq:fxxx} with these values and $z'=\eta$.

%\bibliography{bibliographysingularity}

%\printbibliography
\end{document}